\input amstex
\loadbold
\loadeurm
\loadeusm
\openup2\jot \magnification=1200
\tolerance=1000000

\def\ul{\underline}
\def\volume{{\text{vol}}}
\def\A{\Cal A}
\def\I{\Cal I}
\def\J{\Cal J}
\def\Ok{\Cal O}
\def\d{\partial}
\def\G{\Gamma}
\def\L{\Lambda}
\def\O{\Omega}
\def\cls{{\text{ cls}}}
\def\Z{\Bbb Z}

\def\rank{{\text {rank}}}
\def\mod{{\text{mod\ }}}
\def\norm{{\text{norm}}}
\def\D{\Delta}
\def\H{{\Cal H}}
\def\rad{{\text {rad}}}
\def\e{{\text{e}}}
\def\P{\Cal P}

\def\R{{\Bbb R}}
\def\C{{\Bbb C}}

\def\mult{{\text{mult}}}
\def\diag{\text{diag}}
\def\cls{\text{cls}}
\def\clx{\text{clx}}
\def\K{\Bbb K}
\def\Ok{\Cal O}
\def\M{\Cal M}
\def\a{\alpha}
\def\b{\beta}
\def\g{\gamma}
\def\Q{\Cal Q}
\def\tr{{\text{Tr}}}

\def\Op{{\Cal O}_{\P}}
\def\supp{{\text{supp}}}
\def\ord{{\text{ord}}}
\def\spn{{\text{span}}}

\documentstyle{amsppt}
\pageheight{7.7in}
\vcorrection{-0.05in}
\topmatter
\pageno 1
\title Hecke operators on Hilbert-Siegel
modular forms \endtitle
\author Suzanne Caulk and Lynne H. Walling\endauthor
\address S. Caulk, Department of Mathematics, Regis University,
Denver Colorado \endaddress
\address L.H. Walling, Department of Mathematics,
University of Colorado, Boulder CO 80309\endaddress

\abstract We define Hilbert-Siegel modular forms and Hecke ``operators''
acting on them.  As with Hilbert modular forms (i.e. with Siegel degree 1),
these linear transformations are not linear operators 
until we consider a direct product of spaces
of modular forms (with varying groups), modulo natural identifications
we can make between certain spaces.  With Hilbert-Siegel forms (i.e. with
arbitrary Siegel 
degree) we identify several families of natural identifications
between certain spaces of modular forms.
We associate
the Fourier coefficients
of a form in our product space to even integral lattices, independent
of basis and choice of coefficient rings.  We then 
determine the action of the Hecke operators on these Fourier coefficients,
paralleling the result of Hafner and Walling for Siegel modular forms
(where the number field is the field of rationals).
\endabstract

\endtopmatter
\document
\head{1. Introduction.}\endhead 
A Siegel modular form $F$ of degree $n$
over the rationals has a Fourier series supported
on even integral symmetric
$n\times n$ matrices.  An even integral symmetric matrix can be 
interpreted as the matrix for a quadratic form on an 
even integral lattice,
relative to some $\Z$-basis for that lattice.  Given the transformation
property of $F$ under the symplectic group, the coefficient of $F$
attached to a matrix $T$ is equal to that attached to the conjugate
$^tGTG$ where $G$ is any integral change of basis matrix
(with determinant 1 when $k$, the weight of the modular form,
is odd).  Consequently
we can rewrite $F$ as a ``Fourier series'' supported on 
even integral lattices,
without specifying a basis for each lattice.
For each prime $p$ there are $n+1$ Hecke operators, $T(p)$ and
$T_j(p^2)$ ($1\le j\le n$) associated to $p$, $n$ of which are algebraically
independent.  In [5] we determined the action of these operators on
the Fourier coefficients of $F$.  In this paper we extend this result
to Hilbert-Siegel modular forms.

With $\K$ a totally real number field and $\P$ a prime ideal, we 
mimic the construction of the classical Hecke operators and
construct a linear transformation $T(\P)$ acting on Hilbert modular
forms.  When $\P$ is not principally generated, $T(\P)$ maps modular
forms attached to $\Gamma=SL_2(\Ok)$ ($\Ok$ the ring of integers of $\K$),
to forms attached to the ``psuedo-conjugate''
$$\pmatrix \P&0\\0&1\endpmatrix \Gamma \pmatrix \P^{-1}&0\\0&1\endpmatrix
= \left\{\pmatrix a&b\\c&d\endpmatrix:\ a,d\in\Ok,\ b\in\P,\ c\in\P^{-1},\ 
ad-bc=1\ \right\}.$$
So for $T(\P)$ to be a linear operator
(meaning its domain and codomain are equal), it is necessary to consider
a (finite) direct product of spaces of modular forms attached to 
psuedo-conjuagates of $SL_2(\Ok)$.
In [7], Shimura defined ``Fourier coefficients'' attached to integral
ideals of a form in this direct product, and he determined the action
of $T(\P)$ on these Fourier coefficients.

In the case of Hilbert-Siegel modular forms, we need to consider a (finite)
direct product of spaces of forms attached to psuedo-conjugates of
$Sp_n(\Ok)$ for the maps $T(\P)$ and $T_j(\P^2)$ to be linear operators.
For a form in this direct product, we define ``Fourier coefficients''
attached to even integral lattices, independent of basis and 
choices of coefficient rings 
(note that an $\Ok$-lattice is not necessarily free,
and there are numerous ways to write it as $\A_1 x_1\oplus\cdots\oplus
\A_n x_n$ with the $\A_i$ fractional ideals).  Then we determine the
action of the Hecke operators on these Fourier coefficients.
When $k$ is odd, we need to impose an orientation on $\Lambda$.
Thus 
$$F(\tau)=\sum_{\cls\Lambda} c(\Lambda)e^*\{\Lambda\tau\}$$
where $\cls\L$ runs over isometry classes of lattices $\L$, and
 $e^*\{\Lambda\tau\}=\sum_G \exp\left(\pi i Tr(\ ^tGTG\tau)\right)$; here
$\Lambda=\A_1x_1\oplus\cdots\oplus\A_nx_n$, $T=\big(B(x_i,x_j)\big)$
where $B$ is the symmetric bilinear form 
associated to the quadratic form $Q$
on $\Lambda$ so that $Q(x)=B(x,x)$, 
and $G$ varies
over $O(\Lambda)\backslash GL_n(\Z)$ when $k$ is even, and over
$O^+(\Lambda)\backslash SL_n(\Z)$ when $k$ is odd.  (Two lattices $\L, \O$
are in the same isometry class if there is an isomorphism from one onto the
other that preserves the quadratic form.  Also,
$O(\Lambda)$ is
the orthogonal group of $\Lambda$.)

We begin by defining symplectic groups $\G(\I_1,\ldots,\I_n;\J)$
for fractional ideals $\I_i,\J$.  We show that the spaces of modular
forms associated to $\G(\I_1,\ldots,\I_n;\J)$ and
$\G(\I_1',\ldots,\I_n';\J')$ are naturally isomorphic whenever
$\cls (\I_1\cdots\I_n)=\cls(\I_1'\cdots\I_n')$ and $\clx^+\J=\clx^+\J'$.
(Here $\cls\I$ denotes the wide ideal class of $\I$, and $\clx^+\J$
denotes the strict ideal complex of $\J$.  Thus $\cls\I=\cls\I'$
if $\I=\alpha\I'$ for some $\alpha\in\K$, and $\clx^+\J=\clx^+\J'$
if $\J=\alpha\I^2\J'$ for some fractional ideal $\I$ and $\alpha\gg0$.)
We set $\M_k=\prod_{\I_i,\J} \M_k(\G(\I_1,\ldots,\I_n;\J))/\sim$
(so we identify spaces that are naturally isomorphic).  Next we attach
the Fourier coefficients of (the components of) $F$ to even integral lattices,
independent of the basis and the coefficient rings used to realize each
lattice.  
(For a full discussion of this, see the discussion preceeding Proposition 2.2.)
In \S3
we introduce operators $S(\Q)$ attached to fractional ideals $\Q$,
and we decompose $\M_k$ as $\oplus_{\chi}\M_k(\chi)$ where $\chi$ varies over
ideal class characters, and $F|S(\Q)=\chi(\Q)F$ for $F\in\M_k(\chi)$.
Then in \S4 we
 introduce the Hecke operators $T(\P)$ and $T_j(\P^2)$, $0\le j\le n$, and
we find coset representatives giving the action of the operators.  When then 
analyzing the action of the Hecke operators $T_j(\P^2)$ 
in \S5 we encounter incomplete
character sums; we complete these by replacing $T_j(\P^2)$ with 
$\widetilde T_j(\P^2)$, a combination of $T_\ell(\P^2)$, $0\le \ell\le j$.
Finally, we show that for $\L^\J$ an even integral lattice and 
$F\in\M_k(\chi)$, the $\L^\J$th coefficient of $F|\widetilde T_j(\P^2)$ is
$$\sum_{\P\L\subseteq\O\subseteq\P^{-1}\L} N(\P)^{E_j(\O,\L)}
\chi(\P)^{e_j(\O,\L)} \alpha_j(\O,\L) c_F(\O^\J)$$
where $E_j(\O,\L)$ and $e_j(\O,\L)$ are given by formulas in terms of
the invariant factors $\{\O:\L\}$, and $\alpha_j(\O,\L)$ reflects some
geometry of $(\O\cap\L)/\P(\O+\L)$.  
(A formula for $\alpha_j(\O,\L)$ is given at the end of \S5.)
A similar but much simpler argument
shows that the $\L^\J$th coefficient of $F|T(\P)$ is
$$\sum_{\P\L\subseteq\O\subseteq\L} N(\P)^{E(\O,\L)} c_F(\O^{\J\P^{-1}})$$
(see Theorem 5.2).

In \S 6 we present a lemma on completing a symmetric coprime pair to a
symplectic matrix.
The reader is referred to [6] for basic results on lattices and quadratic
forms.

The authors are thankful for the referee's careful reading.

\head{2. Definitions, isomorphisms, and Fourier coefficients
attached to even integral lattices.}\endhead
Let $\K$ be a totally real number field of degree $d$ over $\Bbb Q$,
and let $\partial$ denote the different of $\K$.  
Let $\H_{(n)}$ denote degree $n$ Siegel upper half-space; so
$$\H_{(n)}=\left\{ \tau=X+iY:\ X,Y\in\R^{n,n} {\text{ are symmetric, }}
Y>0\ \right\}.$$
For fractional ideals $\I_1,\ldots,\I_n,\J$, let
$$\align \G(\I_1,\ldots,\I_n;\J) 
& =\Big\{\pmatrix A&B\\ C&D\endpmatrix\in G_{2n}(\K):\ 
A\ ^tB,C\ ^tD {\text{\ symmetric}},\\
&\qquad\qquad
 A\ ^tD-B\ ^tC=uI,\ 
u\in\Ok^+,\ 
 a_{ij}\in\I_i\I_j^{-1},\\
&\qquad\qquad
 b_{ij}\in\I_i\I_j\J\d^{-1},\ 
c_{ij}\in(\I_i \I_j\J)^{-1}\d,\ d_{ij}\in\I_i^{-1}\I_j\ \Big\}.
\endalign$$
(Here $A=\pmatrix a_{ij}\endpmatrix$, etc.)
So with $\G=\G(\Ok,\ldots,\Ok;\Ok)$, $\G(\I_1,\ldots,\I_n;\J)$ corresponds
to the formal conjugate $\ul\delta \G\ul\delta^{-1}$ where
$$\ul\delta=
\pmatrix
\J\I_1\\&\ddots\\&&\J\I_n\\&&&
\I_1^{-1}\\&&&&\ddots\\&&&&&\I_n^{-1}\endpmatrix.$$
Also notice that $\G(\I_1,\ldots,\I_n;\J\I^2)=\G(\I_1\I,\ldots,\I_n\I;\J)$.

\smallskip
\noindent{\bf Definition.}  A degree $n$ ($n>1$), weight $k$ Hilbert-Siegel
modular form for $\G(\I_1,\ldots,\I_n;\J)$ is a function
$f:\H^d_{(n)}\to\C$ so that the following two conditions hold.
\itemitem{(1)} $f$ is ``analytic on $\H^d_{(n)}$ and at infinity,'' meaning
that for $\tau\in\H^d_{(n)}$,
$$f(\tau)=\sum_{T} c(T) e\{T\tau\}$$
where $T$ runs over symmetric, positive semi-definite
$n\times n$ matrices. 
Also,
$\sigma(M)$ denotes the trace of a matrix $M$,
$\tr$ denotes the trace from $\K$ to $\Q$, and
$$e\{T\tau\}=\exp(\pi i\sigma(\tr(T\tau))).$$
Here $\tr(T\tau)=\sum_{i=1}^d T^{(i)}\tau_i$, 
where $T^{(i)}$ is the image of $T$
under the $i$th embedding of $\K$ into $\R$.
\itemitem{(2)}  For all $M\in\G(\I_1,\ldots,\I_n;\J)$,
$f|M=f$
where, for any matrix $\pmatrix A&B\\C&D\endpmatrix$ (written in $n\times n$
blocks), we define
$$\align
&f|\pmatrix A&B\\C&D\endpmatrix (\tau)\\
& =\det(N(A\ ^tD-B\ ^tC))^{k/2} \det(N(C\tau+D))^{-k}\\
&\qquad
f((A\tau+B)(C\tau+D)^{-1}).\endalign$$
\noindent Here $N$ denotes the norm from $\K$ to $\Bbb Q$, extended so that
$$N(C\tau+D)=\prod_{i=1}^d C^{(i)}\tau_i+D^{(i)}.$$

Let $\M_k(\G(\I_1,\ldots,\I_n;\J))$ denote
the space of Hilbert-Siegel modular forms for $\G(\I_1,\ldots,\I_n;\J)$,
and let $f$ be a modular form in this space.
Since $f(\tau+B)=f(\tau)$ for all symmetric $B\in (\I_i\I_j\J\partial^{-1})$,
$\partial$ the different of $\K$,
we must have $e\{TB\}=1$ for all $T\in\supp f$.  Note that for
$T=(t_{ij})$, $B=(b_{ij})$ symmetric matrices,
$$\sigma(TB)=\sum_{i=1}^n t_{ii} b_{ii} \ + \ 
\sum_{1\le i<j\le n} 2 t_{ij}b_{ij}.$$
Thus for $T\in\supp f$, we must have $T\in \left ( (I_iI_j\J)^{-1}
\right)$ with $T$ even, meaning $t_{ii}\in 2 I_i^{-2}\J^{-1}.$

\smallskip\noindent
{\bf Definitions.}  We define families of isomorphisms between spaces
of modular forms as follows:
Fix $f\in\M_k(\G(\I_1,\ldots,\I_n;\J))$.

First, for $\a\in\K^{\times}$ and $1\le\ell\le n$, let
$$M=\pmatrix I_{\ell-1}\\&\a^{-1}\\&&I_{n-\ell}\endpmatrix,$$
and  define
$$f|U_{\ell}(\a)= f|\pmatrix M\\&^tM^{-1} \endpmatrix.$$
Since 
$\pmatrix M\\&^tM^{-1}\endpmatrix
\G(\I_1',\ldots,\I_n';\J)
\pmatrix M^{-1}\\&^tM\endpmatrix
=\G(\I_1,\ldots,\I_n;\J)$
where $$\I_i'=\cases \I_i&\text{if $\ell\not= i$}\\
\a\I_{\ell}&\text{if $\ell=i$,}\endcases$$
$U_{\ell}(\a)$ defines an isomorphism from
$\M_k(\G(\I_1,\ldots,\I_n;\J))$ onto $\M_k(\G(\I_1',\ldots,\I_n';\J))$.

For $\a\gg0$, define $W(\a):\M_k(\G(\I_1,\ldots,\I_n;\J))\to
\M_k(\G(\I_1,\ldots,\I_n;\a\J))$ by
$$f|W(\a)= f|\pmatrix \a^{-1}I_n\\&I_n\endpmatrix.$$
One easily checks (as we did above for $U_\ell(\alpha)$) that $W(\alpha)$
is an isomorphism.

For $\Q$ a fractional ideal, $1\le \ell<n$, choose
$$A\in\pmatrix \Q^{-1}&\Q\I_{\ell}\I_{\ell+1}^{-1}\\
\Q^{-1}\I_{\ell}^{-1}\I_{\ell+1}&\Q \endpmatrix$$
so that $\det A=1$ (possible by Strong Approximation; see p. 42 [6]).
Let
$$M=\pmatrix I_{{\ell}-1}\\&A\\&&I_{n-\ell-1}\endpmatrix,$$
and define
$$f|V_{\ell}(\Q)=f|\pmatrix M\\& ^tM^{-1} \endpmatrix.$$
Since $\pmatrix M\\&^tM^{-1}\endpmatrix
\G(\I_1',\ldots,\I_n';\J)
\pmatrix M^{-1}\\&^tM\endpmatrix
=\G(\I_1,\ldots,\I_n;\J)$
where $$\I_i'=\cases \I_i&\text{if $i\not= \ell,\ell+1$}\\
\Q\I_{\ell}&\text{if $i=\ell$}\\
\Q^{-1}\I_{\ell+1}&\text{if $i=\ell+1$,}\endcases$$
the map $V_{\ell}(\a)$ defines an isomorphism from
$\M_k(\G(\I_1,\ldots,\I_n;\J))$ onto $\M_k(\G(\I_1',\ldots,\I_n';\J))$.

For $\Q$ a fractional ideal and $1\le \ell<j\le n$, set
$$V_{\ell j}(\Q)=V_\ell(\Q)V_{\ell+1}(\Q)\cdots V_{j-1}(\Q).$$
Then 
$V_{\ell j}(\Q)$ defines an isomorphism from
$\M_k(\I_1,\ldots,\I_n;\J))$ onto $\M_k(\G(\I_1',\ldots,\I_n';\J))$, where
$$\I_i'=\cases \I_i&\text{if $i\not= \ell,j$}\\
\Q\I_{\ell}&\text{if $i=\ell$}\\
\Q^{-1}\I_{j+1}&\text{if $i=j$,}\endcases.$$

\proclaim{Proposition 2.1} The maps $U_{\ell}(\a)$, $W(\a)$, 
$V_{\ell j}(\Q)$ commute and for fixed $\ell, j$, these operators are 
multiplicative (as functions on fractional ideals). 
\endproclaim

\demo{Proof}  The tedious aspect of proving such relations among our
isomorphisms is that, for any of the above listed maps, the domain and
codomain differ.  Keeping track of appropriate domains and codomains, and
using the matrices we used to define the actions of these operators,
it is then straightforward to verify the operators commute, remembering
that if $MN^{-1}\in\G$ for any group
$\G=\G(\I_1,\ldots,\I_n;\J)$, then $f|M=f|N$ for $f\in\M_k(\G)$.
$\square$
\enddemo

\noindent{\bf Definition.} For $f\in\M_k(\G(\I_1,\ldots,\I_n;\J))$,
$g\in\M_k(\I_1',\ldots,\I_n';\J'))$,
define the equivalence relation $\sim$ by
$f\sim g$  if some composition of the maps
$U_i, W, V_{ij}$ takes $f$ to $g$.  We define
$$\M_k=\prod_{\I_i,\J} \M_k(\G(\I_1,\ldots,\I_n;\J))/\sim$$
where $\I_1,\ldots,\I_n,\J$ vary over all fractional ideals.  Note that
$\sim$ partitions the spaces $\M_k(\G(\I_1,\ldots,\I_n;\J))$ according
to $\cls (\I_1\cdots\I_n)$, $\clx^+\J$.
Thus $\M_k\approx \prod_{\cls\I,\clx^+\J} \M_k(\G(\Ok,\ldots,\Ok,\I;\J))$, 
$\cls\I$ runs over all ideal classes
and $\cls^+\J$ runs over all strict ideal class complexes.
($\J,\J'$ are in the same strict ideal class complex if $\J'=\a\I^2\J$ for
some fractional ideal $\I$ and $\a\gg0$.)
\smallskip

Given an element $F\in\M_k$, we can associate the Fourier coefficients of
$F$ with lattices equipped
with positive semi-definite, 
even integral quadratic forms as described below.

First note the following.
Say $f\in\M_k(\G(\I_1,\ldots,\I_n;\J))$
is a component of a chosen representative for
$F$.  So the support of $f$ lies in
$\big((\I_i\I_j\J)^{-1}\big)$.
For $T\in\supp f$, consider $T$ as defining a quadratic form $Q$
on $\Ok x_1\oplus\cdots\oplus\Ok x_n$ (recall that $T$ is symmetric
and even, meaning the $i$th diagonal entry of $T$ lies in 
$2\I_i^{-2}\J^{-1}$).  
Thus with $\L=\I_1 x_1\oplus\cdots\oplus\I_n x_n$, $Q(\L)\subseteq 2\J^{-1}$.
If $\alpha\Ok=\J$ then the scaled lattice $\L^{\alpha}$ is even integral,
meaning $\alpha Q(\L)\subseteq 2\Ok$.  When $\J$ is not principal, we
abuse notation and language
 and refer to $\L^{\J}$ as an even integral lattice.
We agree to identify
 $\Lambda^{\I^2}$ with $\I\Lambda$ since locally
everywhere these are identified (paralleling the fact that
$\G(\I_1,\ldots,\I_n;\I^2\J) = \G(\I\I_1,\ldots,\I\I_n;\J)$).

Not every lattice is a free $\Ok$-module, and given a lattice $\L$,
there are many ways to choose fractional ideals $\I_j$ and vectors $x_j$
so that $\L=\I_1 x_1\oplus\cdots\oplus\I_n x_n$.
Say we have 
$$\L=\I_1 x_1\oplus\cdots\oplus\I_n x_n
=\I'_1 y_1\oplus\cdots\oplus\I'_n y_n,$$
$y_j=\sum_i a_{ij}x_i$.  Then by 81:8 [6],
$\I_1\cdots\I_n=\I'_1\cdots\I'_n\cdot\det(a_{ij})$.

The Invariant Factor Theorem (81:11 [6]) says that given lattices
$\L, \O$ on a (non-zero) space $V$, there are vectors $x_1,\ldots,x_n\in V$
and fractional ideals $\I_1,\ldots,\I_n,\A_1,\ldots,\A_n$ so that
$$\align
\L&=\I_1 x_1\oplus\cdots\oplus\I_n x_n,\\
\O&=\I_1\A_1 x_1\oplus\cdots\oplus\I_n\A_n x_n,
\endalign$$
and $\A_i|\A_{i+1}$ ($1\le i<n$);
the $\A_i$ are unique and called the invariant factors of $\O$ in $\L$.
We use $\{\L:\O\}$ to refer to these invariant factors, and we write
$\{\L:\O\}=(\A_1,\ldots,\A_n)$.

When analyzing the action of Hecke operators on Fourier coefficients,
we sum over lattices $\L$ where $\P\L\subseteq\O\subseteq\P^{-1}\L$,
$\P$ a prime ideal and $\L$ a fixed reference lattice of rank $n$.
By the Invariant Factor Theorem, we have sublattices $\L_i$ so that
$$\align
\L&=\L_0\oplus\L_1\oplus\L_2,\\
\O&=\P\L_0\oplus\L_1\oplus\P^{-1}\L_2.
\endalign$$
So for instance, $r_0=\rank\L_0$ is the multiplicity of $\P$ among the
invariant factors $\{\L:\O\}$, denoted $r_0=\mult_{\{\L:\O\}}(\P)$.

We will also need to consider 
$(\L^{\J}\cap\O^{\J})/\P(\L^{\J}+\O^{\J})\approx \L_1^{\J}/\P\L_1^{\J}.$
We will only be considering even integral $\L^{\J}$.  
Thus $Q$ induces a quadratic form ${1\over 2}\alpha Q$ on 
$\L_1^{\J}/\P\L_1^{\J}$
defined by 
$${1\over 2}\alpha Q(x+\P\L_1)={1\over 2}\alpha Q(x)+\P \in\Ok/\P$$
where $\alpha\in\K$ has been fixed so that $\alpha\Ok_{\P}=\J\Ok_{\P}$.
Since $Q(\L_1)\subseteq 2\J^{-1}$, this gives us a
quadratic form on the $\Ok/\P$-space $\L_1^{\J}/\P\L_1^{\J}$.
Note that the structure of the quadratic space 
$\L_1^{\J}/\P\L_1^{\J}$ is independent of the choice of $\alpha$.

\smallskip
\noindent{\bf Definition.}  Let $f\in\M_k$.
Given any even integral positive semi-definite 
lattice $\L^{\J}$ with
$\L=\I_1x_1\oplus\cdots\oplus\I_nx_n$, we set
$$c(\L^{\J})=c_F(\L^{\J})=c_f(T)\cdot N(\I_1\cdots\I_n)^k N(\J)^{nk/2}  $$
where $f\in\M_k(\G(\I_1,\ldots,\I_n;\J))$ is a representative of the
component of $F$ corresponding to $\cls(\I_1\cdots\I_n)$,
$\clx^+\J$, and
$T=\pmatrix B(x_i,x_j)\endpmatrix.$  
If $k$ is odd, we assume $\L$ is also equipped with an
orientation.
\smallskip

\proclaim{Proposition 2.2}  The ``Fourier coefficient''
$c(\L^{\J})$ is well-defined.
\endproclaim

\demo{Proof}
First, suppose we also have $\L'=\I_1 y_1\oplus\cdots\oplus\I_n y_n$,
and $\L=\L'$.  Take $M=(\alpha_{ij})$ so that
$$(y_1\ldots y_n)=(x_1\ldots x_n)M.$$
Hence $(B(y_i,y_j))=\ ^tMTM$ (recall $T=(B(x_i,x_j))$).
Note that 
$$\I_j y_j = \sum_i \alpha_{ij}\I_j x_i \subseteq \L,$$
so $\alpha_{ij}\in\I_i \I_j^{-1}$. 
 Also, since $\volume \L=\volume\L'$, it
follows that $\det M\in\Ok^{\times}$.  (Recall that if $k$ is odd then
$\L$ has an orientation, and $\det M$ must also be totally positive.)
Thus $\pmatrix M\\&^tM^{-1}\endpmatrix \in \G(\I_1,\ldots,\I_n;\J)$,
and so
$$f=f|\pmatrix M\\&^tM^{-1}\endpmatrix.$$
Hence
$$\align
f(\tau)&=\sum_T c_f(T) e\{T\tau\}\\
&= N(\det M)^k \sum_T c_f(T) e\{TM\tau \ ^tM\}\\
&=\sum_T c_f(T) e\{\ ^tMTM\tau\}
\endalign$$
(recall that $\det M$ is a unit, and that if $k$ is odd, a totally
positive unit).
Thus $c_f(\ ^tMTM)=N(\det M)^k\ c_f(T)=c_f(T)$, and so
$$\align
c(\L^{\J})& = c_f(T) N(\I_1\cdots\I_n)^k N(\J)^{nk/2}\\
& = c_f(\ ^tMTM)  N(\I_1\cdots\I_n)^k N(\J)^{nk/2}
= c(\L^{\prime \J}).
\endalign$$
Thus the definition of $c(\L^{\J})$ is independent of the choice
of basis relative to the coefficient ideals $\I_1,\ldots,\I_n$ and
the scaling ideal $\J$ (so here $\I_1,\ldots,\I_n$ and $\J$ are fixed).

Next, fix $\J$ and suppose
$\L=\I_1 x_1\oplus\cdots\oplus\I_n x_n
=\I'_1 y_1'\oplus\cdots\oplus\I'_n y_n'.$
Then by 81:8 of [6], $\I'_1\cdots\I'_n\in\cls (\I_1\cdots\I_n).$
Thus, as we have seen, 
$M_k(\G(\I_1,\ldots,\I_n;\J))\simeq
M_k(\G(\I_1',\ldots,\I_n';\J))$ via an appropriate composition of
the maps $U_i,V_{ij}$; the action of this composition is given by
$$f\mapsto f'=  f|\pmatrix M\\&^tM^{-1}\endpmatrix$$
where $M\in\left(\I_i\I_j^{-1}\right)$ and
$(\det M)\I'_1\cdots\I'_n=\I_1\cdots\I_n$.
Also,
$$f'(\tau)= N(\det M)^k f(M\tau\ ^tM)
= N(\det M)^k \sum_T c_f(T) \e\{\ ^tMTM\tau\},$$
so $c_{f'}(^tMTM)= N(\det M)^k c_f(T)$.
Set $(y_1\ldots y_n)=(x_1\ldots x_n)M;$ thus
$\Ok y_1\oplus\cdots\oplus\Ok y_n\simeq \ ^tMTM.$
We claim that with $\L'=\I_1'y_1\oplus\cdots\oplus\I_n'y_n$, $\L'=\L$.
We know that writing $M$ as $(\alpha_{ij})$,
$$y_j=\sum_{i=1}^n \alpha_{ij}x_i \in(\I_j')^{-1}(\I_1x_1\oplus\cdots
\oplus\I_nx_n)$$
and so $\L'\subseteq\L$.  Since 
$\norm(\L)=(\I_1\cdots\I_n)^2\cdot\det T=(\I_1'\cdots\I_n')^2\cdot
\det(^tMTM)=\norm(\L'),$
we have $\L'=\L$.  
So
$$\align
c(\L^{\prime\J})& = c_{f'}(\ ^tMTM) \ N(\I'_1\cdots\I'_n)^k\ N(\J)^{nk/2}\\
&=c_f(T)\ N(\det M)^k\ N(\I'_1\cdots\I'_n)^k\ N(\J)^{nk/2}\\
&=c_f(T)\ N(\I_1\cdots\I_n)^k\ N(\J)^{nk/2}\\
&=c(\L^{\J}).
\endalign$$
Thus by our assumption, $\I_1'y_1\oplus\cdots\oplus\I_n'y_n =
\I_1'y_1'\oplus\cdots\oplus\I_n'y_n',$
which reduces the problem to the preceeding case.  Hence our definition
of $c(\L^{\J})$ is independent of the choice of classes $\I_1,\ldots,\I_n$
so that $\cls(\I_1\cdots\I_n)$ is as prescribed (so here $\cls(\I_1\cdots\I_n),
\J$ are fixed).

Finally, suppose $\J'\in\clx^+\J$.  Thus $\J'=\alpha\I^2\J$ for some
fractional ideal $\I$ and some $\alpha\gg0$.  Say
$\L=\I_1 x_1\oplus \cdots \oplus \I_n x_n.$
We have agreed previously to identify $\L^{\I^2}$ and $\I\L$, so
that $c_f(\L^{\alpha\I^2\J})=c_f(\I\L^{\alpha\J}),$
whether we think of $f$ as associated to $\G(\I_1,\ldots,\I_n;\alpha\I^2\J)$
or to $\G(\I\I_1,\ldots,\I\I_n;\alpha\J)$ (remember, these are two names
for the same group).
So suppose $\J'=\alpha\J$,
$\alpha\gg0$.  
We know $\M_k(\G(\I_1,\ldots,\I_n;\J))\simeq 
\M_k(\G(\I_1,\ldots,\I_n;\alpha\J))$ via
$$f\mapsto f'=f|\pmatrix \alpha^{-1}I\\&I\endpmatrix.$$
Thus
$$f'(\tau) = N(\alpha)^{-nk/2}\ f(\alpha^{-1}\tau)
= N(\alpha)^{-nk/2}\ \sum_T c_f(T) e\{\alpha^{-1} T\tau\}.$$
Consequently $c_{f'}(\alpha^{-1}T)=N(\alpha)^{-nk/2}\ c_f(T)$.
Set $\L'=\I_1 x_1\oplus\cdots\oplus\I_n x_n$ equipped with the quadratic
form $\alpha^{-1}T$ (so $(\L')^{\alpha\J}$ is an integral lattice).
Then we have
$$\align
c(\L^{\J}) &= c_f(T)\ N(\I_1\cdots\I_n)^k\ N(\J)^{nk/2} \\
&=c_{f'}(\alpha^{-1}T)\ N(\I_1\cdots\I_n)^k\ N(\alpha\J)^{nk/2} \\
&= c((\L')^{\alpha\J}).
\endalign$$
Thus the definition of $c(\L^{\J})$ is also independent of the choice
of the representative for $\clx^+\J$.
$\square$
\enddemo

\head{3. Eigenspaces.}\endhead
The operators $V_{\ell}(\Q)$ are lifts of the 
Hilbert modular form operators that Eichler called
$V(\Q^{-1})$ [3] and Shimura called $S(\Q)$ [7].  We now introduce
another lift of these operators, and following
Shimura [7] (where
$n=1$), we decompose $\M_k$ into eigenspaces for these new operators.

\smallskip\noindent
{\bf Definition.}
Let $\Q$ be a fractional ideal, and choose 
$$\pmatrix a&b\\ c&d \endpmatrix
\in \pmatrix \Q&\Q^{-1}\I_\ell^2\J\d^{-1}
\\ \Q\I_\ell^{-2}\J^{-1}\d&\Q^{-1} \endpmatrix$$
with $ad-bc=1$.
Set $$M=\pmatrix I_{\ell-1}&&&0_{\ell-1}\\&a&&&b\\&&I_{n-\ell}&&&0\\
0_{\ell-1}&&&I_{\ell-1}\\&c&&&d\\&&0&&&I_{n-\ell} \endpmatrix.$$
Then $M\G(\I_1,\ldots,\I_n;\J)M^{-1}=\G(\I_1',\ldots,\I_n';\J)$ where
$$\I_j'=\cases I_j&\text{if $j\not=\ell$,}\\
\Q^{-1}\I_\ell&\text{if $j=\ell.$} \endcases$$
Thus $S_{\ell}(\Q):\M_k(\G(\I_1,\ldots,\I_n;\J))\to
\M_k(\G(\I_1',\ldots,\I_n';\J))$ is an isomorphism where we define
$$f|S_{\ell}(\Q)=f|M.$$

\proclaim{Proposition 3.1}
With $\Q,\P$ fractional ideals and $\a\in\K^{\times}$,  $S_\ell(\Q)$
commutes with $U_i(\a)$, $W(\a),$ and  $V_{ij}(\P)$. 
Further, $S_i(\Q)V_{ij}(\Q)=S_j(\Q)$ and
$U_i(\alpha^{-1}) S_i(\Q)=S_i(\alpha\Q)$.
\endproclaim

\demo{Proof}
Keeping in mind the various domains for different incarnations of our
functions, it's easy to show $S_{\ell}(\Q)$ commutes with $U_i(\a), W(\a)$.

To show $S_\ell(\Q),V_{ij}(\P)$ commute, it suffices to show $S_\ell(\Q)$,
$V_i(\P)$ commute (recall how $V_{ij}$ is defined).  When $\ell\not=i,i+1$,
it is easy to see the matrices giving the actions of $S_\ell(\Q), V_i(\P)$ 
commute.  As we explain below, we can reduce our attention to the case
$n=2$.

Suppose $\ell=i$ or $i+1$; for the sake of clarity, let us first look
at the case where $i=1, \ell=1$ or 2.  Then the action of each operator
is given by a matrix of the form
$$\pmatrix A&&B\\&I&&0\\C&&D\\&0&&I \endpmatrix$$
where $A,B,C,D$ are $2\times 2$ matrices.  Since the product of such
matrices (and their inverses) is again of this form, it suffices to
restrict our attention to the submatrices $\pmatrix A&B\\C&D\endpmatrix.$

First consider $n=2$, $i=\ell=1$.  Choose
$$\pmatrix a&b\\c&d\endpmatrix \in
\pmatrix \Q&\Q\P^2\I_1^2\J\d^{-1}\\ (\Q\I_1^2\J)^{-1}\d&\Q^{-1} \endpmatrix$$
so that $ad-bc=1$.  Choose
$$A\in\pmatrix \P^{-1}&\P\Q\I_1\I_2^{-1}\\
\P^{-1}\I_1^{-1}\I_2&\P \endpmatrix $$
so that $\det A=1$.  (Note that these choices are possible, even when
$\Q=\P$.)  Then
$$M=\pmatrix a&&b\\&1&&0\\c&&d\\&0&&1 \endpmatrix $$
gives the action of both
$$S_1(\Q):\M_k(\G(\Q\I_1,\I_2;\J))\to\M_k(\G(\I_1,\I_2;\J))$$
and
$$S_1(\Q):\M_k(\G(\P\Q\I_1,\P^{-1}\I_2;\J))
\to\M_k(\G(\P\I_1,\P^{-1}\I_2;\J)).$$
Similarly,
$N=\pmatrix A\\&^tA^{-1}\endpmatrix $
gives the action of both
$$V_1(\Q):\M_k(\G(\Q\I_1,\I_2;\J))\to\M_k(\G(\P\Q\I_1,\P^{-1}\I_2;\J))$$
and
$$V_1(\Q):\M_k(\G(\I_1,\I_2;\J))
\to\M_k(\G(\P\I_1,\P^{-1}\I_2;\J)).$$

A simple (but tiresome) check shows $MNM^{-1}N^{-1}\in\G(\Q\I_1,\I_2;\J)$,
which implies that $S_1(\Q),V_1(\P)$ commute when $n=2$.  For general $n$,
we have $S_1(\Q)V_1(\P)S_1(\Q^{-1})V_1(\P^{-1})$ represented by a matrix
of the form 
$$M'=\pmatrix A&&B\\&I&&0\\C&&D\\&0&&I\endpmatrix$$
where $\pmatrix A&B\\C&D\endpmatrix\in\G(\Q\I_1,\I_2;\J)$.  Thus
$M'\in\G(\Q\I_1,\I_2,\ldots,\I_n;\J)$, and $S_1(\Q), V_1(\P)$ commute
for general $n$.  Similarly, $S_2(\Q),V_1(\P)$ commute for $n=2$, and thus
for general $n$.

For general $n,\ell,i$ with $\ell=i$ or $i+1$, the matrices giving the
action of $S_\ell(\Q),V_i(\P)$ are of the form
$$\pmatrix I_{i-1}&&&0_{i-1}\\ &A&&&B\\ &&I&&&0\\
0_{i-1}&&&I_{i-1}\\&C&&&D\\&&0&&&I \endpmatrix,$$
where $A,B,C,D$ are $2\times 2$
matrices.  So again the problem reduces to showing 
$S_\ell(\Q),V_i(\P)$ ($\ell=i$ or $i+1$)
 commute when $n=2$ (which we have done).

To see that $S_i(\Q)V_{ij}(\Q)=S_j(\Q)$, it again suffices to consider
$j=i+1$. First, fix a group $\G=\G(\I_1,\ldots,\I_n;\J)$.
From our definitions, we have a product of three matrices giving
the action of
$S_i(\Q)V_i(\Q)S_{i+1}(\Q^{-1})$ on $\M_k(\G).$
One easily verifies that the
conditions on these matrices ensure the product for
$S_i(\Q)V_i(\Q)S_{i+1}(\Q^{-1})$ lies in $\G$.  
Finally, one verifies $U_i(\alpha^{-1})S_i(\Q)=S_i(\alpha\Q)$
by matrix multiplication.
$\square$
\enddemo

This shows that the $S_\ell(\Q)$ act on $\M_k$, where, for $F\in\M_k$,
$F \sim (\ldots,f_i,\ldots)$, 
$F|S_\ell(\Q) \sim (\ldots,f_i|S_\ell(\Q),\ldots).$
It also 
shows that $S_i(\Q), S_j(\Q)$
are equivalent on $\M_k$,
and so we simply refer to this operator on $\M_k$ as $S(\Q)$.  
Furthermore, on $\M_k$, $S(\Q)=S(\a\Q)$
for all $\a\in\K^{\times}$; since we also know $S(*)$ is multiplicative,
the map $\cls\I\mapsto S(\I)$ gives a group action of
the ideal class group  on $\M_k$.

\proclaim{Proposition 3.2}  $\M_k=\displaystyle{\oplus_{\chi}\M_k(\chi)}$
where $\chi$ varies over all characters of the ideal class group, and
$$\M_k(\chi)=\{F\in\M_k:\ F|S(\Q)=\chi(\Q) F \text{ for all }\Q\ \}.$$
\endproclaim

\demo{Proof}  First notice that $\M_k(\chi)\cap\M_k(\psi)=\{0\}$ if
$\chi\not=\psi$.  To prove this take $F\in\M_k(\chi)\cap\M_k(\psi)$.
Then $\chi(\Q)F=F|S(\Q)=\psi(\Q)F$ for all $\Q$.  
Since $\chi\not=\psi$ there is a $\Q$
such that $\chi(\Q)\not=\psi(\Q)$.  Therefore $F=0.$

For $\psi$ an ideal class character, let 
$\displaystyle{G_\psi={1\over h}\sum
_{\cls\I}\overline{\psi}(\I)\ F|S(\I)}$ 
where $h$ is the class number of $\K$.  Note that
$$\sum_{\psi}G_\psi={1\over h} \sum_{\cls\I}
\left(\sum_{\psi}\overline\psi(\I)\right) F|S(\I)=F
\ \text{since} \ 
\sum_{\psi}\overline\psi(\I)=\cases h&\text{if $\cls\I=\cls\Ok$,}\\
0&\text{otherwise}.\endcases$$  
Thus $F=\sum_{\psi}G_\psi.$

Next notice that
$$\align
G_\psi|S(\Q)
&={1\over h}
\sum_{\cls\I}\overline\psi(\I)\ F|S(\I)|S(\Q)\\
&={1\over h}
\sum_{\cls\I}\overline\psi(\I)\ F|S(\I\Q)\\
&={1\over h}
\psi(\Q)\sum_{\cls\I}\overline\psi(\I\Q)\ F|S(\I\Q)\\
&=\psi(\Q)\ G_\psi.\endalign$$
Thus $G_\psi\in\M_k(\psi)$, and hence $F\in\oplus_{\chi}\M_k(\chi)$.
$\square$
\enddemo

\head{4. Hecke operators.}\endhead
We begin by defining the Hecke operators.  Then we show that they act on
each $\M_k(\chi)$.  After this, we describe how to find a set of coset
representatives giving the action of the operators.  Finally, 
in the next section we analyze
the action of the operators on Fourier coefficients attached to
even integral lattices, proving our main
theorem.

\smallskip\noindent
{\bf Definition.}
Let $\P$ a prime ideal; set $\G=\G(\I_1,\ldots,\I_n;\J)$
and $\G'=\G(\I_1,\ldots,\I_n;\P\J).$
We define the Hecke operator
$T(\P): \M_k(\G')\to \M_k(\G)$ 
by
$$F|T(\P)=N(\P)^{n(k-n-1)/2}\sum_{\gamma} F|\gamma$$
where $\gamma$ runs over a complete set of coset representatives for
$(\G'\cap\G)\backslash\G$.  
Note that $\G'$ is the formal conjugate of $\G$ by the matrix
${\delta}=\pmatrix \P I_n\\& I_n\endpmatrix$.  When $\K=\Bbb Q$, we
define $T(p)$ on $\M_k(\G)$ by
$$f|T(p)=p^{n(k-n-1)/2}\sum_{\gamma} f|\delta^{-1}\gamma$$
where $\gamma$ runs over a complete set of coset representatives for
$(\G'\cap\G)\backslash\G$, $\G'=\delta\G\delta^{-1}$,
$\delta=\pmatrix pI\\&I\endpmatrix$.
(This normalization of $T(p)$ is standard, and as with the standard 
normalization of the degree 1 Hecke operator $T(p)$, the purpose of the
normalization is to force the coefficient to be 1 on the ``lead'' term
in the expression for the $\L$th coefficient of $F|T(p)$.)

Now fix $1\le j\le n$; let 
$\G'_j=\G(\P\I_1,\ldots,\P\I_j,\I_{j+1},\ldots,\I_n;\J).$
We define
the Hecke operators $T_j(\P^2):\M_k(\G'_j)\to\M_k(\G)$ by
$$F|T_j(\P^2)= 
\sum_{\gamma}F|\gamma$$
where $\gamma$ runs over a complete set of
coset representatives for $(\Gamma'_j\cap\Gamma)\backslash\Gamma$. 

Note that $\G'_j$ is the formal conjugate of $\G$ by
$\diag(\P I_j,I_{n-j},\P^{-1}I_j,I_{n-j})$.
When $\K=\Bbb Q$, $\P=p\Z$,
we define $T_j(\P^2)=T_j(p^2)$ on $\M_k(\G)$ by
$$f|T_j(\P^2)= \sum_{\gamma} f|\delta^{-1}\gamma$$
where $\delta=\diag(pI_j,I_{n-j},{1\over p}I_j,I_{n-j})$,
$\G'=\delta\G\delta^{-1}$, and $\gamma$ runs over a complete set of
coset representatives for $(\G'\cap\G)\backslash\G$.
(We introduce a normalization later.)

\proclaim{Proposition 4.1} The operators $T(\P), T_j(\P^2)$ 
commute with $U_i(\a), W(\b), V_{i\ell}(\Q)$, and $S_i(\Q)$ where
$\a,\b\in\K^{\times}$ with $\b\gg0$, and $\Q$ is a fractional ideal.
Thus $T(\P)$, $T_j(\P^2)$ act on $\M_k(\chi)$ (as defined in
Proposition 3.2).
\endproclaim

\demo{Proof}
To show $T_j(\P^2)$ commutes with the $V_{i\ell}(\Q)$, it suffices
to show it commutes with $V_i(\Q)$.  Take
$$A\in\pmatrix \Q^{-1}&\I_i\I_j^{-1}\Q\P\\
\I_i^{-1}\I_j\Q^{-1}\P&\Q \endpmatrix$$
with $\det A=1$, and set
$$\I_i'=\cases \P\I_i&\text{if $i\le j$,}\\
\I_i&\text{otherwise.}\endcases$$
Then with 
$$M=\pmatrix I_{i-1}\\&A\\&&I_{n-i-1}\\&&&I_{i-1}\\&&&&^tA^{-1}\\
&&&&&I_{n-i-1} \endpmatrix,$$
$M$ gives the action of $V_i(\Q):\M_k(\G)\to\M_k(M^{-1}\G M)$ and of
$V_i(\Q):\M_k(\G')\to\M_k(M^{-1}\G' \M)$
(note that these are the appropriate codomains).  
Now let $\{\gamma\}$ be a complete set of coset representatives
for $(\G'\cap\G)\backslash\G$.  Thus $\{M^{-1}\gamma M\}$ is a complete
set of coset representatives for 
$(M^{-1}\G'M\cap M^{-1}\G M)\backslash M^{-1}\G M$.  Hence for
$f\in\M_k(\G')$,
$$\align
f|T_j(\P^2)V_i(\Q)
&= \sum_{\gamma} f|\gamma|M \\
&= \sum_{\gamma} f|M|M^{-1}\gamma M \\
&= f|V_i(\Q)T_j(\P^2).\endalign$$

Similarly, to show $T_j(\P^2)$ commutes with $S_i(\Q)$, choose
$$\pmatrix a&b\\c&d\endpmatrix \in
\pmatrix \Q&\Q^{-1}\I_i^2\P^2\J\d^{-1}\\
\Q\I_i^{-2}\J^{-1}\d&\Q^{-1} \endpmatrix$$
so that $ad-bc=1$.  Then $\pmatrix a&b\\c&d\endpmatrix$ lifts to a matrix
$M$ so that $M$ gives the action of $S_i(\Q):\M_k(\G)\to
\M_k(M^{-1}\G M)$ and of $S_i(\Q):\M_k(\G')\to\M_k(M^{-1}\G' M)$.
Thus 
$$\align
f|T_j(\P^2)S_i(\Q)
&= \sum_{\gamma} f|\gamma|M \\
&= \sum_{\gamma} f|M^{-1}|M^{-1}\gamma M \\
&= f|S_i(\Q)T_j(\P^2).\endalign$$
Similar but simpler arguments show that $T_j(\P^2)$ commutes with
$U_i(\a), W(\b)$, and that $T(\P)$ commutes with $V_i(\Q),
S_i(\Q), U_i(\alpha), W(\beta)$.  $\square$
\enddemo

\proclaim{Proposition 4.2}  For $f'\in\M_k(\G'),$ we have
$$f'|T_j(\P^2)
=  \sum_{\O,\L_1,Y} f'|S^{(j)}(\O)
\pmatrix I&Y\\&I\endpmatrix \pmatrix C^{-1}\\&^tC\endpmatrix.$$
Here $\O$ varies over all lattices such that
$\P\L\subseteq\O\subseteq\P^{-1}\L$, $\overline\L_1$ varies over all
codimension $n-j$ subspaces of $\O\cap\L/\P(\O+\L)$, $C=C(\O,\L_1)$.
With $r_0=\mult_{\{\L:\O\}}(\P)$, $m_1=\mult_{\{\L:\O\}}(\Ok)$,
$r_1=m_1-n+j$, $\mu\in\P^{-1}\ -\ \O$ fixed,
$$S^{(j)}(\O)=\left(\prod_{i=r_0+1}^{r_0+r_1}S_i(\P)\right)
\left(\prod_{i=r_0+r_1+1}^{j}S_i(\P^2)\right),$$
$$Y=\pmatrix W_0&W_2&0&W_3\\^tW_2&\mu W_1\\0\\^tW_3\endpmatrix$$
with $Y\in(\I_i\I_j\J\partial^{-1})$,
$W_0$ varying modulo $\P^2$, $W_1,W_2,W_3$ varying modulo $\P$ with
$\P$ not dividing $\det W_1$.
(Recall that $\partial$ is the different of $\K$.)
  Here $W_0$ is $r_0\times r_0$ and
symmetric, $W_1$ is $r_1\times r_1$ and symmetric, $W_2$ is
$r_0\times r_1$, $W_3$ is $r_0\times(n-j)$.
\endproclaim

\demo{Proof}
Let $\P$ be a prime ideal, and
fix $j$, $1\le j\le n$.
We essentially follow the algorithm presented in [5] to find a set of
coset representatives giving the action of
$$T_j(\P^2):\M_k(\G(\P\I_1,\ldots,\P\I_j,\I_{j+1},\ldots,\I_n;\J))
\to \M_k(\G(\I_1,\ldots,\I_n;\J)).$$
For convenience, we will take $\I_i=\Ok$ for $1\le i\le j$,
$\I_i=\Ok$ for $j<i<n$; also, we take $\I=\I_n$ and 
$\J\partial^{-1}$ to be integral ideals relatively prime to $\P$
(recall that the equivalence class of
$\M_k(\G)$ is determined by $\cls\I,\cls^+\J$, and $\partial$ is the
different of $\K$).  Note that this allows us to choose $\mu$
relatively prime to $\I_1\cdots\I_n\J\partial^{-1}=\I\J\partial^{-1}$.

Choose $M\in\G(\Ok,\ldots,\Ok,\I;\J)$, and let $M_j=(A|B)$ denote the
top $j$ rows of $M$ with $A, B$ $j\times n$ matrices.
Let $\L=\Op x_1\oplus \cdots \oplus \Op x_n$ be a reference lattice.

{\bf Step 1.}  Let 
$$\O_0=\ker(\L\to \L(A)\ \mod \P\Op)$$
where $A=(a_1\cdots a_n)$ and $\L\to \L(A)\ \mod \P\Op$ denotes the map
that takes $x_i$ to $\overline a_i$ (which is a $1\times j$ matrix with entries
in $\Op/\P\Op$).  Note that the $r_0=\rank_{\P}\L(A)$ is at most $j$
since $A$ is a $j\times n$ matrix.

We claim there is a matrix $\pmatrix C\\&^tC^{-1}\endpmatrix \in
\G(\Ok,\ldots,\Ok,\I;\J)$ such that 
$\O_0=\L C_0\pmatrix \P  I_{r_0}\\&I_{n-r_0}\endpmatrix$
and $(A|B)\pmatrix C\\&^t C^{-1}\endpmatrix
=(A'|B')$ with $a_1',\ldots,a_n'\equiv 0\ (\mod \P)$.

First, write $A=(a_1\ldots a_n)$ and consider the rank modulo $\P\Op$ of
$(a_1\ldots a_n)$.  Let $E_1$ be an $(n-1)\times(n-1)$ invertible matrix
(i.e. a change of basis matrix) so that
$$(a_1\ldots a_n)E_1=(a'_1\ldots a'_n)$$
with $a'_1,\ldots,a'_{r_0}$ linearly independent modulo $\P\Op$ and
$a'_{r_0+1}\equiv\cdots\equiv a'_{n-1}\equiv 0\ (\mod \P\Op).$
Note that 
$G_1=\pmatrix E_1\\&^tE_1^{-1}\\&&1\endpmatrix\in
\G(\Ok,\ldots,\Ok,\I;\J).$
If $a_n$ is in the span modulo $\P\Op$ of $a'_1,\ldots,a'_{r_0}$, then there
is a matrix $E_2=\pmatrix I&*\\&1\endpmatrix$ such that
$(a'_1\ldots a'_{n-1} a_n)E_2=(a'_1\ldots a'_n)$
where $a'_n\equiv 0\ (\mod \P\Op);$ note that 
$G_2=\pmatrix E_2\\&^t E_2^{-1}\endpmatrix \in \G(\Ok,\ldots,\Ok,\I;\J)$
and we take $C_0$ to be $G_1 G_2$.
If $r_0=n-1$ then we are now done, regardless of whether $a_n$ is in the
span (modulo $\P\Op$) of $a'_1,\ldots,a'_{r_0}$.

So suppose $r_0<n-1$ (and thus $a'_{n-1}\equiv 0\ (\mod \P\Op)$) and
$a_n$ is not in the span (modulo $\P\Op$) of $a'_1,\ldots,a'_{r_0-1}$.

Choose $\eta\in\I^{-1}, \rho\in\I$ such that $\eta\equiv\rho\equiv 1\ (\mod\P)$
and choose $\nu\in\P$ such that $(\nu,\eta\rho)=1$.  Thus there are 
$\a,\b\in\Ok$ so that $\a\nu-\rho\eta\b=1$.  Then with 
$E_3=\pmatrix I\\&\matrix \a&\eta\b\\ \rho&\nu\endmatrix \endpmatrix$,
$(a'_1\ldots a'_{n-1} a_n)E_3=(a'_1\ldots a'_{n-2} a''_{n-1} a''_n )$
with $a''_{n-1}\equiv a_n\ (\mod\P)$, $a''_n\equiv0\ (\mod\P).$
Note that $G_3=\pmatrix E_3\\&^tE_3^{-1}\endpmatrix\in
\G(\Ok,\ldots,\Ok,\I;\J).$
Let $E_4$ be the permutation matrix that permutes columns $r_0$ and $n-1$;
then $G_4=\pmatrix E_4\\&^tE_4^{-1}\endpmatrix \in
\G(\Ok,\ldots,\Ok,\I;\J)$
and $(a_1\ldots a_n)G_1 G_3 G_4 \equiv (a'_1\ldots a'_{r_0-1} a_n 0 \ldots 0)
\ (\mod \P).$
Hence in this case we take $C_0=G_1 G_3 G_4.$

Thus there is a matrix $C_0$ and integer $r_0$ such that
$\pmatrix C_0\\&^tC_0^{-1}\endpmatrix \in \G(\Ok,\ldots,\Ok,\I;\J)$
and $\O_0=\L C_0\pmatrix \P  I_{r_0}\\&I_{n-r_0}\endpmatrix.$
Then with renewed notation, $M_j
\pmatrix C_0\\&^tC_0^{-1}\endpmatrix$ has the form
$(a_1,\ldots,a_n|b_1,\ldots,b_n)=(A_0 A_1|B)$, $A_1\equiv 0\ (\mod \P)$.

Note that while $C_0$ is not uniquely determined, $\O_0$ is.

{\bf Step 2.}  First note that Lemma 7.2 of [5] easily generalizes to
number fields, where we ``permute'' $b_\ell$ and $b_n$ as we ``permuted''
$a'_{r_0}$ and $a_n$ in the preceeding paragraph.  Thus with
$$M_j\pmatrix C_0\\&^tC_0^{-1}\endpmatrix = (A|B)
=(a_1,\ldots,a_n|b_1,\ldots,b_n),$$
where $b_1,\ldots,b_{r_0}$ are in the span mod $\P\Op$ of $a_1,\ldots,a_{r_0}$,
and the rank mod $\P\Op$ of $(a_1,\ldots,a_{r_0},b_{r_0+1},\ldots,b_n)$
is $j$.  Thus for some $C=\pmatrix I_{r_0}\\&E\endpmatrix$ with
$\pmatrix C\\&^tC^{-1}\endpmatrix\in \G(\Ok,\ldots,\Ok,\I;\J)$, we have
$(A|B)\pmatrix C\\&^tC^{-1}\endpmatrix=(A'|B')$ with $a_i'=a_i$ for
$i\le r_0$, and $j$ the rank mod $\P\Op$ of $(a_1,\ldots,a_{r_0},
b'_{r_0+1},\ldots,b'_j)$.
We want to accomplish the above rearrangement, as well as replacing
$b'_{j+1},\ldots,b_n'$ with vectors in the span mod $\P\Op$ of
$a_1,\ldots,a_{r_0}$; we want to identify these modifications with a 
uniquely determined lattice.

From Step 1 we have $\L=\L_0\oplus\D_1$, $\O_0=\P\L_0\oplus\D_1$ with
$\D_1$ uniquely determined modulo $\P\L$.  This corresponds to a splitting
$\L^\#=\L_0'\oplus\L_1'$ of the (formal) dual of $\L$, where $\L_0'$ is 
orthogonal to $\D_1$ and $\L_1'$ is orthogonal to $\L_0$.  
(So $$\L^{\#}=\Op y_1\oplus \cdots \oplus \Op y_n$$
and the basis $\{y_1,\ldots,y_n\}$ is dual to $\{x_1,\ldots,x_n\}$.)

Let $V$ be the $\Op/\P\Op$-space consisting of all $j\times 1$ matrices,
and let $U$ 
be the subspace spanned by $\overline a_1,\ldots,\overline a_{r_0}$.
Let
$$\O_1'=\ker(\L^{\#}\to \L^{\#}(B)\ \mod \P\Op \to V/U),$$
where $\L^{\#}\to \L^{\#}(B)$ corresponds to $y_i\mapsto b_i$,
and the map into $V/U$ is the canonical projection map.  Thus as in
Step 1, we can find a matrix $C=\pmatrix I_{r_0}\\&E\endpmatrix$
so that $\pmatrix C\\& ^t C^{-1}\endpmatrix \in \G(\Ok,\ldots,\Ok,\I;\J)$
and $(A|B)\pmatrix C\\& ^t C^{-1}\endpmatrix =(A'|B')$ with
$b_0'=b_0,\ldots,b_{r_0}'=b_{r_0},b'_{j+1},b'_n$ in the span mod $\P\Op$
of $a_0,\ldots,a_{r_0}$.  Also, setting $C_1=C_0 C$,
$$\O'_1= \L^{\#}\ ^tC_1^{-1} 
\pmatrix I_{r_0}\\ &\P I_{j-r_0}\\ && I_{n-j} \endpmatrix ;$$
$\L^{\#}=\L^{\#}\ ^t C_1^{-1} = \L_0'\oplus \L_2'\oplus \L_3'$ with
$\rank \L_2'=j-r_0$ and $\L_0'\oplus\L'_3$ uniquely determined modulo
$\P\L^{\#}$.
Correspondingly,
$$\O_0=\L C_1\pmatrix \P I_{r_0}\\& I_{n-r_0}\endpmatrix.$$
Note that $\L=\L C_1=\L_0\oplus\D_2\oplus\L_3$.  Set
$$\O_1=\L C_1\pmatrix \P I_{r_0}\\& I_{j-r_0}\\&& \P I_{n-j}\endpmatrix
=\P\L_0\oplus\D_2\oplus\P\L_3.$$
Since $\L_0'\oplus\L_3'$ are uniquely determined modulo $\P\L^{\#}$, $\D_2$
is uniquely determined modulo $\P\L$.

{\bf Step 3.}  Write $$M_j\pmatrix C_1\\& ^tC_1^{-1}\endpmatrix
=(a_1\cdots a_n|b_1\cdots b_n) = (A_0 A_1 A_3| B_1 B_2 B_3)$$
where $j$ is the rank modulo $\P\Op$ of $(A_0,B_1)$,
$A_1,A_3\equiv 0\ (\mod \P\Op)$, and $B_0,B_3$ are in the (column) span
modulo $\P\Op$ of $A_0$.  We want to modify $A_1$ to be of the form
$(A_1',A_2')$ where $A_2'\equiv 0\ (\mod \P^2\Op)$.  Recall that we have
$$\L=\L_0\oplus\D_2\oplus\L_3, \qquad
\O_1=\P\L_0\oplus \D_2\oplus \P\L_3,$$
with $\rank \L_0=r_0, \rank\D_2=j-r_0$.  Renewing our notation, let
$(x_1,\ldots,x_n)$ be a basis corresponding to this decompositon of $\L$.

Recall we have fixed $\mu\in\P^{-1}-\Ok$ so that $\mu$ is relatively
prime to $\I\J\partial^{-1}$; set
$$\P\O_2=\ker(\O_1\to \O_1(\mu A)\ \mod \P\Op)$$
where $\O_1\to \O_1(\mu A)$ denotes the map taking $x_i$ to
$\mu a_i$, so $\mu\O_1(A)\ \mod \P\Op$ is spanned by
$\overline a_1,\ldots,\overline a_{r_0},
\mu \overline a_{r_0+1},\ldots,\mu \overline a_j$
with $\overline a_1,\ldots,\overline a_{r_0}$ linearly independent modulo
$\P\Op$ (recall $a_i\simeq 0\ (\mod \P\Op)$ for $i>j$).
Thus 
$$\P\O_2=\P^2\L_0\oplus\P\L_1\oplus\L_2\oplus\P\L_3$$
where $\L_2$ is uniquely determined modulo $\P\O_1$.
As in the previous steps, we can find a matrix
$C=\pmatrix I_{r_0}\\&E\\&& I_{n-j}\endpmatrix$ such that
$\pmatrix C\\& ^tC^{-1}\endpmatrix \in \G(\Ok,\ldots,\Ok,\I;\J)$ and
$$\P\O_2=\L C \pmatrix \P^2 I_{r_0}\\&\P I_{r_1}\\&&I_{r_2}\\&&&\P I_{n-j}
\endpmatrix.$$
Correspondingly, $(A|B)\pmatrix C\\& ^tC^{-1}\endpmatrix 
=(A_0 A_1' A_2 A_3 | B_0 B_1'  B_2 B_3)$
with $A_1'\equiv 0\ (\mod \P\Op)$, $A_2\equiv 0\ (\mod \P^2\Op)$, and
$a_1,\ldots,a_{r_0}, \mu a'_{r_0+1},\ldots,\mu a'_{r_0+r_1}$
linearly independent modulo $\P\Op$ where $r_1=\rank \L_1$ ( and so
$A_1'$ is $j\times r_1$).  Let $C(\O,\L_1)=C_0 C_1.$

{\bf Step 4.}
Write $M_j\pmatrix C_2\\&^tC_2^{-1}\endpmatrix = 
(A_0,A_1,A_2,A_3|B_0,B_1,B_2,B_3).$
So $A_1,A_3\equiv 0\ (\mod \P\Op)$, $A_2\equiv 0\ (\mod \P^2\Op)$, and
the columns of $(A_0,\mu A_1)$ are linearly independent modulo
$\P\Op$.  Also, $B_0,B_3$ are in the column span modulo $\P\Op$ of $A_0$.

Since $B_0, B_1$ are in $\spn_{\P}A_0$, we can solve
$$A_0 Y_0'\equiv -B_0\ (\mod \P\Ok_{\P}),
\qquad A_0 Y_3\equiv -B_3\ (\mod \P\Ok_{\P}).$$
Note that as $B\ ^tA$ is symmetric and 
$A_1,A_2,A_3\equiv 0 \ (\mod \P\Ok_{\P})$,
$B_0\ ^tA_0$ is symmetric modulo $\P\Op$.  Also, since $A_0$ has full rank
modulo $\P$, there is some matrix $E\in GL_j(\Op)$ such that
$E A_0=\pmatrix A'\\0\endpmatrix.$  Writing $E B_0=\pmatrix B'\\B''\endpmatrix,$
we see $B''\equiv 0\ (\mod \P\Op)$ since $E(B_0\ ^tA_0)\ ^tE$ is symmetric
modulo $\P\Op$.  Thus $Y_0'$ is the unique solution modulo $\P$ to 
$A' Y_0'\equiv -B'\ (\mod \P\Op);$ since
$(A')^{-1}\ B'$ is symmetric modulo $\P\Op$, we can choose $Y_0'$ to
be symmetric.

Let 
$$\align
&(A_0,A_1,A_2,A_3|B_0',B_1,B_2,B_3')\\
&\qquad =
(A_0,A_1,A_2,A_3|B_0,B_1,B_2,B_3)
\pmatrix I&&&&Y_0&0&0&Y_3\\&I&&&0\\
&&I&&0\\
&&&I&^tY_3\\
&&&&I\\&&&&&I\\&&&&&&I\\
&&&&&&&I\endpmatrix,
\endalign$$
$B_0',B_3'\equiv0\ (\mod \P\Ok_{\P})$.
Then just as we argued about $Y_0'$, there is a unique modulo $\P\Op$
symmetric solution $Y'$ to
$$(A_0,\mu A_1) Y' \equiv -(\mu( B_0'+A_3\ ^tY_3),B_1)\ (\mod \P\Op).$$
Decompose $Y'$ as $\pmatrix Y_0''&Y_2\\^tY_2&Y_1\endpmatrix;$
choose $\delta\in\P$ so that $\delta\mu\equiv 1\ (\mod \P)$
and
set $Y_0=Y_0'+\delta Y_0''$.

Note that since $\rank_{\P}(A_0,B_1)=r_0+r_1$, we have
$\rank_{\P}(B_1+A_0 Y_0)=\rank_{\P} B_1=r_1.$
Since $-\mu A_1 Y_1\equiv B_1+A_0 Y_2\ (\mod \P\Ok_{\P})$, we must have
$\det Y_1\in\Op^{\times}$.

Take $Y=\pmatrix W&0&W_3\\0&I\\^tW_3\endpmatrix
=\pmatrix W_0&W_2&0&W_3\\^tW_2&\mu W_1\\0&&I\\^tW_3\endpmatrix$
to be a symmetric matrix in $\left(\I_i\I_j\J\partial^{-1}\right )$ with
$W_0\equiv Y_0\ (\mod \P^2\Ok_{\P})$ and
$W_i\equiv Y_i\ (\mod \P\Ok_{\P})$ for $i=1,2,3$.

Then 
$$(A_0,A_1,A_2,A_3|B_0,B_1,B_2,B_3)
\pmatrix I&Y\\0&I\endpmatrix
=(A_0,A_1,A_2,A_3|B_0'',B_1',B_2',B_3')$$
with $B_0''\equiv 0\ (\mod \P^2\Ok_{\P})$, 
$B_2'\equiv B_2\ (\mod\P^2\Ok_{\P})$,
$B_1', B_3'\equiv 0\ (\mod\P\Ok_{\P})$.

Let $C=C(\O,\L_1)$.  Also,
identifying $S_i(\P)$ with a matrix giving its action, let
$$S^{(j)}(\O)=\left(\prod_{i=r_0+1}^{r_0+r_1}S_i(\P)\right)
\left(\prod_{i=r_0+r_1+1}^{j}S_i(\P^2)\right).$$
We see that, with renewed notation,
$$M_j\pmatrix C\\&^tC^{-1}\endpmatrix
\pmatrix I&Y\\&I\endpmatrix S^{-1}(\O)
=(A_0,A_1,A_2,A_3|B_0,B_1,B_2,B_3)$$
with $A_3,B_3\equiv 0\ (\mod\P\Op)$, $B_0,B_1,B_2\equiv0\ (\P^2\Op)$.

However, while the matrices for $S_i(\P^{-2})$ lie in $\G$, the matrices
for $S_i(\P^{-1})$ do not.  We remedy this as follows.

For $r_0<i\le r_0+r_1$, choose $\a_i\in\P^{-1}$, 
$\b_i\in\P\I_i^2\J\partial^{-1}$, $\g_i\in\P^{-1}\I_i^{-2}\J^{-1}\partial$
so that $\a_i\delta-\b_i\g_i=1$ and for any prime $\Q\not=\P$ dividing
$\delta$, $\Q$ does not divide $\b_i$.
(Recall that our choice of $\delta$ ensures $\delta\in\P-\P^2$.)  
Set $\underline\a=\diag(\ldots,\a_i,\ldots)$, an $r_1\times r_1$
matrix; define $\underline\b,
\underline\g$ in an analogous fashion.
So
$$\pmatrix I_{r_0}&&&0_{r_0}\\
&\underline\a&&&\underline\b\\
&&I&&&0\\ 
0_{r_0}&&&I_{r_0}\\&\underline\g&&&\delta I\\
&&0&&&I\endpmatrix$$
gives the action of $\prod_{i=r_0+1}^{r_0+r_1} S_i(\P^{-1}).$

Now consider 
$$\pmatrix I&-\mu W_1\\&I\endpmatrix
\pmatrix \underline\a&\underline\b\\ \underline\g&\delta I\endpmatrix
\equiv\pmatrix \underline\a-\mu W_1\underline\g&\underline\b- W_1\\
\underline\g&\delta I\endpmatrix\ (\mod\P).$$
We find that $\pmatrix \underline \b- W_1\\ \delta I\endpmatrix$
is a coprime symmetric right-hand pair for 
$\G(\I_{r_0+1},\ldots,\I_{r_0+r_1'};\J)$
($2r_1'\times 2r_1'$ matrices).  Thus by Lemma 6.1,
there exist matrices $U,V$ so that
$$\pmatrix U&\underline\b+\delta W_1\\ V&\delta I\endpmatrix
\in \G(\I_{r_0+1},\ldots,\I_{r_0+r_1'};\J).$$
Hence
$$X^{-1}=
\pmatrix I_{r_0}&&&0_{r_0}\\
&U&&&\underline\b- W_1\\
&&I&&&0\\
0_{r_0}&&&I_{r_0}\\
& V&&&\delta I\\&&0&&&I
\endpmatrix
\in\G,$$
and with $Y'=\pmatrix W_1&W_2&W_3'\\^tW_2\\^tW_3'\endpmatrix$ and
$$N^{-1}=\pmatrix C\\&^tC^{-1}\endpmatrix
\pmatrix I&-Y'\\
&I \endpmatrix
X^{-1}\left(\prod_{i=r_0+r_1'+1}^j S_i(\P^{-2})\right),$$
we have $M_j N^{-1}=(A_0,A_1,A_2,A_3|B_0,B_1,B_2,B_3)$ with
$A_3,B_3\equiv0\ (\mod\P\Op)$, $B_0,B_1,B_2\equiv0\ (\mod\P^2\Op)$.
Thus by an easy generalization of Lemma 7.1 of [5],
$MN^{-1}\in\G\cap\G'$.  Also, since
$$\pmatrix \delta I&-\underline\b+ W_1\\-\underline\g&\underline\a
\endpmatrix
\pmatrix U&\underline\b-W_1\\V&\delta I\endpmatrix
=\pmatrix I&0\\V'&I\endpmatrix
\in\G(\I'_{r_0+1},\ldots,\I'_{r_0+r_1'};\J),$$
we have
$$f'|N=f'|S(\O)\pmatrix I&-Y\\&I\endpmatrix \pmatrix C^{-1}\\&^tC\endpmatrix.
\ \square$$
\enddemo

\proclaim{Proposition 4.3}  For $f'\in\M_k(\G(\I_1,\ldots,\I_n;\P\J))$,
we have
$$f'|T(\P)=N(\P)^{n(k-n-1)/2} \sum_{\O,Y_0} f'| S(\O)
\pmatrix I&&Y_0\\&I&&0\\&&I\\&&&I\endpmatrix
\pmatrix C^{-1}\\&^tC\endpmatrix$$
and $Y_0\in(\I_i\I_{\ell}\J\partial^{-1})$ varies over symmetric
$r\times r$ matrices modulo $\P$ and $C$ varies as in Proposition 4.2.
Here $r=\mult_{\{\L:\O\}}(\P)$ and $S(\O)=\prod_{i=r+1}^n S_i(\P)$.
\endproclaim

\demo{Proof}
To find coset representatives for $T(\P)$, take $M\in\G(\I_1,\ldots,\I_n;\J)$;
write $M=\pmatrix A&B\\C&D\endpmatrix.$  Let
$\O=\ker(\L\mapsto\L(A)\ (\mod\P)),$ and choose $C=C(\O)$ so that
$\O=\L C\pmatrix \P I_r\\&I_{n-r}\endpmatrix.$  Thus
$(A|B)C=(A_0 A_1|B_0 B_1)$, $A_1\equiv 0\ (\mod\P\Op)$, $\rank_\P A_0=r$
where $A_0$ is $n\times r$, and $B_0\in\spn_\P A_0$.  Choose symmetric
$Y_0$ with $i,\ell$-entry in $\I_i\I_\ell\J\partial^{-1}$ such that
$A_0 Y_0\equiv B_0\ (\mod\P\Op)$.  Then with 
$Y=\pmatrix Y_0\\&\mu I\endpmatrix$,
$(A|B)\pmatrix C\\&^t C^{-1}\endpmatrix \pmatrix I&-Y\\&I\endpmatrix
S^{-1}(\O) = (A'|B')$
where $S(\O)=\prod_{i=r+1}^n S_i(\P)$ and $B'\equiv 0\ (\mod\P\Op)$.
As before, choose diagonal $(n-r)\times(n-r)$ matrices $\underline\a$,
$\underline\b$, $\underline\g$ so that
the action of $S(\O)$ is given by $\pmatrix I_r&&0_r\\&\underline\a
&&\underline\b\\0_r&&I_r\\&\underline\g&&\delta I\endpmatrix.$
Then
$$\pmatrix I&-Y\\&I\endpmatrix \pmatrix I&&0\\&\underline\a&&\underline\b\\
0&&I\\&\underline\g&&\delta I\endpmatrix
=\pmatrix I&&-Y\\&I&&0\\&&I\\&&&I\endpmatrix
\pmatrix I&&0\\&\underline\a+\mu\underline\g&&\underline\b+\delta I\\
0&&I\\&\underline\g&&\delta I\endpmatrix.$$
Here $\underline\b\equiv 0\ (\mod\P\Op)$, $\rank_\P(\underline\b+\mu\delta I)
=n-r$.  Thus $\pmatrix \underline\b-\mu\delta I\\ \delta I\endpmatrix$
is a symmetric coprime right-hand pair for $\G(\I_{r+1},\ldots,\I_n;\J)$,
hence there are $U, V$ so that $\pmatrix I&&0\\&U&&\underline\b+\mu\delta I\\
0&&I\\&V&&\delta I\endpmatrix \in \G(\I_1,\ldots,\I_n;\J)$.
Since $\pmatrix U&\underline\b+\mu\delta I\\V&\delta I\endpmatrix
\pmatrix \delta I&-\underline\b\\-\underline\g&\underline\a\endpmatrix
\in \G(\I_{r+1},\ldots,\I_n;\P\J)$ we get the result as claimed.
$\square$
\enddemo

\head{5. Evaluating the action of the Hecke operators.}\endhead

When evaluating the action of the operators $T_j(\P^2)$, we encounter
incomplete character sums.  To complete these, we define modified operators
as follows.

\smallskip\noindent
{\bf Definition.}  For $\P$ a prime ideal and $1\le j\le n$, define
$$\widetilde T_j(\P^2)=N(\P)^{j(k-n-1)}\sum_{0\le\ell\le j}
\beta(n-\ell,j-\ell)
S_{\ell+1}(\P)\cdots S_j(\P) T_\ell(\P^2).$$

We will also need the following rather technical result.

\proclaim{Proposition 5.1}  
Let $T_1$ be a symmetric $r_1\times r_1$ matrix whose $i,\ell$-entry lies in
$(\I_{r_0+i}\I_{r_0+\ell}\J)^{-1}$, and whose $i$th diagonal entry
lies in $2\I_{r_0+i}^{-2}\J^{-1}$;
fix $\mu\in \P^{-1} - \Ok$.
With $W$ varying
over all symmetric $r_1\times r_1$ matrices modulo $\P$ with $i,\ell$-entry
in $\I_{r_0+i}\I_{r_0+\ell}\J\partial^{-1}$,
$$\sum_W \e\{\mu T_1 W\} = \sum_{0\le m\le r_1'} \sum_{\D,U}
\e\{\mu T_{\D} U\}$$
where for each $m$, $\overline\D$ varies over dimension $m$ subspaces
of $\overline\L_1$, $\overline\D\simeq T_\D\ (\mod\P)$, and
$U$ varies over all $m\times m$ symmetric matrices modulo $\P$.
\endproclaim

\demo{Proof}
For a moment, let's fix $W$.
Since $W$ is symmetric, we can view it as the matrix of a quadratic form on an 
$r_1'$ dimensional $\Ok/\P$ space $V=L/\P L$, 
$L=\left(\I_{r_0+1}^{-1}y_1\oplus\cdots\oplus\I_{r_0+r_1}y_{r_1}\right)
^{\J^{-1}\partial}.$
(When $\P$ is dyadic, let $W$ define an integral quadratic form
on 
$$L=(\Ok_{\P}\I^{-1}_{r_0+1}y_1\oplus
\cdots\oplus\Ok_{\P}\I^{-1}_{r_0+r_1}y_{r_1})
^{\J^{-1}\partial},$$
and let $V=L/\P L$, a quadratic space over $\Ok_{\P}/\P\Ok_{\P}\approx
\Ok/\P$.  We use \S93 of [6] to understand the structure of $L$ and
thereby of $V$.)
The radical of this space is uniquely
defined, so for some $G\in GL_{r_1}(\Ok)$,
$$^t G^{-1} W G^{-1} \equiv \pmatrix U\\&0\endpmatrix\ (\mod \P)$$
where $U$ is $m\times m$ with $\rank_{\P}U=m$.  (So $V G^{-1}=J\oplus
\rad V$ where $J$ is a regular space whose isometry class is uniquely 
determined by $V$, and $J\simeq U$.)

So 
$$\align
\e\{\a T_1 W\}&= \e\{\a T_1\ ^tG\pmatrix U\\&0\endpmatrix G\} \\
&= \e\{\a(G T_1\ ^tG)\pmatrix U\\&0\endpmatrix \} \\
&= \e\{\a SU\}
\endalign$$
where $G T_1\ ^tG=\pmatrix S&*\\*&*\endpmatrix$, $S$ an $m\times m$ matrix.
Here we take $\L$ to be a rank $n$ lattice as in the previous section,
and we equip $\L$ with a quadratic form such that $\L\simeq T$; thus
with $\L_1$ as in the previous section, we have
$\L_1\simeq T_1$, and $\D=\L_1 \ ^tG\pmatrix I_m\\0\endpmatrix \simeq S.$
So $S$ corresponds to an $m$-dimensional subspace $\overline\D$ of the
$\Ok/\P$-space $\overline\L_1$.  Thus each $W$ gives rise to (at least one)
pair $(\overline \D,U)$, $\overline \D$ an $m$-dimenstional subspace of 
$\overline\L_1$, $U$ an $m\times m$ integral symmetric matrix of rank $m$ 
modulo $\P$.  

With $T_1$ still
fixed, fix $m$, $0\le m\le r_1$.  We now define a map 
$\varphi$ from all pairs $(\overline \D,U)$ as above to symmetric
$r_1\times r_1$ matrices $W$.  Here $\overline \D$ is an $m$-dimensional
subspace of $\overline \L_1$, and $U$ is an integral symmetric $m\times m$
matrix with $\rank_{\P}U=m$.  For each such $\overline\D$ we fix some
$G=G_{\D}\in GL_{r_1}(\Ok)$ so that $\overline\D=\overline\L_1\ ^tG
\pmatrix I_m\\0\endpmatrix$.  We define $\varphi(\overline\D,U)=
\ ^tG\pmatrix U\\&0\endpmatrix G.$

We first show that the image of $\varphi$ consists of all symmetric
$r_1\times r_1$ matrices $W$ modulo $\P$ with $\rank_{\P}W=m$.
Then we show that $\varphi$ is injective.

As shown above, given any $W$ in the codomain of $\varphi$,
$$W\equiv \ ^tG\pmatrix U\\&0\endpmatrix G\ (\mod \P)$$
where $U$ is $m\times m$, $m=\rank_{\P} U$, and
$G\in GL_{r_1}(\Ok)$.

Take $\overline \D=\overline \L_1\ ^tG \pmatrix I_m\\0\endpmatrix.$
So $\overline \D$ is an $m$-dimensional subspace of $\overline\L_1$,
and thus $^tG\pmatrix I_m\\0\endpmatrix$ and 
$^tG_{\D}\pmatrix I_m\\0\endpmatrix$
each map a basis for $\overline \L_1$ to a basis for $\overline\D$.
Hence with $(x_1,\ldots,x_{r_1})$ a basis for $\overline\L_1$,
$(x_1,\ldots,x_{r_1})\ ^tG=  (y_1,\ldots,y_{r_1})$,
$(x_1,\ldots,x_{r_1})\ ^tG_{\D}=  (z_1,\ldots,z_{r_1})$,
we must have $(y_1,\ldots,y_m)\equiv (z_1,\ldots,z_m)\ ^tC\ (\mod \P)$
for some $C\in GL_m(\Ok)$.  Thus
$(x_1,\ldots,x_{r_1})\ ^tG = (x_1,\ldots,x_{r_1})\ ^tG_{\D}
\pmatrix \ ^tC&*\\0&*\endpmatrix$, meaning
$^tG=\ ^tG_{\D}\pmatrix \ ^tC&*\\0&*\endpmatrix$.

Hence modulo $\P$, 
$$W\equiv \ ^tG\pmatrix U\\&0\endpmatrix G
=\ ^tG_{\D} \pmatrix \ ^tC&*\\0&*\endpmatrix
\pmatrix U\\&0\endpmatrix
\pmatrix C&0\\*&*\endpmatrix G_{\D}
=\varphi(\overline\D,\ ^tCUC).$$
Thus $\varphi$ is surjective.

Now we show $\varphi$ is injective.  Say
$$W\equiv \varphi(\D_1,U_1)\equiv \varphi(\D_2,U_2)\ (\mod \P).$$
Thus with $G_i=G_{\D_i}$, we have
$$W\equiv \ ^tG_1\pmatrix U_1\\&0\endpmatrix G_1
\equiv \ ^tG_2\pmatrix U_2\\&0\endpmatrix G_2 \ (\mod \P).$$
So with $G=G_2G_1^{-1}$,
$$\pmatrix U_1\\&0\endpmatrix \equiv \ ^tG\pmatrix U_2\\&0\endpmatrix G
\ (\mod \P).$$
Since the columns of $U_i$ are linearly independent modulo $\P$, we must have
$G\equiv \pmatrix C&0\\*&*\endpmatrix \ (\mod \P).$

Now we compare $\D_1,\D_2$; we will find that $\overline\D_1=\overline\D_2$,
so $G_1=G_2$ and hence $U_1\equiv U_2\ (\mod \P).$
With notation as before, we have
$$\align
(y_1,\ldots,y_m) & =(x_1,\ldots,x_{r_1'})\ ^tG_1 \pmatrix I_m\\0\endpmatrix,\\
(z_1,\ldots,z_m) & =(x_1,\ldots,x_{r_1'})\ ^tG_2 \pmatrix I_m\\0\endpmatrix\\
& =(x_1,\ldots,x_{r_1'})\ ^tG_1\ ^tG \pmatrix I_m\\0\endpmatrix.
\endalign$$
Thus given our knowledge of $G$, we see that
$(y_1,\ldots,y_m)=(z_1,\ldots,z_m)\ ^tC$, and hence
$\overline \D_1=\overline \D_2$.
Thus $G_1=G_2$, and consequently $U_1\equiv U_2\ (\mod \P).$
Therefore $\varphi$ is injective. $\square$
\enddemo

We can now prove our main result.
In the remark following the proof we demonstrate how to compute the
geometric term $\alpha_j(\O,\L)$.

\proclaim{Theorem 5.2}  Let $F\in\M_k(\chi)$ where $\chi$ is a character
of the ideal class group and $\M_k(\chi)$ is as defined in Proposition 3.2. 
\item{(1)}  The $\L^\J$th coefficient of
$F|\widetilde T_j(\P^2)$ is
$$\sum_{\P\L\subseteq\O\subseteq\P^{-1}\L}
N(\P)^{E_j(\L,\O)} \chi(\P)^{e_j(\L,\O)} \alpha_j(\O,\L) c_F(\O^\J)$$
where $E_j(\L,\O)=k(r_2-r_0+j)+r_0(r_0+m_1+1)+r_1(r_1+1)/2-j(n+1)$,
$e_j(\L,\O)=2r_2+r_1=r_2-r_0+j$, and $\alpha_j(\O,\L)$ is the number of totally
isotropic codimension $n-j$ subspaces of $\O\cap\L/\P(\O+\L)$.
Here $r_0=\mult_{\{\L:\O\}}(\P)$,
 $m_1=\mult_{\{\L:\O\}}(\Ok)$, $r_1=m_1-n+j$, and 
$r_2=\mult_{\{\L:\O\}}(\P^{-1})$. 
\item{(2)}  The $\L^\J$th coefficient of $F|T(\P)$ is
$$\sum_{\P\L\subseteq\O\subseteq\L} N(\P)^{E(\O,\L)} \chi(\P)^{n-r}
c_F(\O^{\J\P^{-1}})$$
where  $r=\mult_{\{\L:\O\}}(\Ok)$ and $E(\O,\L)=k(n-r)+r(r+1)/2-n(n+1)/2$.
\endproclaim

\demo{Proof}
Take fractional ideals $\I_1,\ldots,\I_n,\J$ and
$f'\in\M_k(\G(\I_1',\ldots,\I_n';\J))$,
$$\I_i'=\cases \P\I_i&\text{if $i\le j$,}\\ \I_i&\text{if $i>j$.}
\endcases$$
In the preceeding proposition, consider the subsum where we
fix a choice of $\O$:
$$\sum_{\L_1,Y} f'| S(\O)\pmatrix I&Y\\&I\endpmatrix 
\pmatrix C^{-1}\\&^tC\endpmatrix
=\chi(\P)^{2r_2+r_1} \sum_{\L_1,Y} f''|\pmatrix I&Y\\&I\endpmatrix 
\pmatrix C^{-1}\\&^tC\endpmatrix$$
where $f''$ is the component of $F$ corresponding to the group
$\G(\I''_1,\ldots,\I_n'';\J)$,
$$\I_i''=\cases \P\I_i&\text{if $1\le i\le r_0$,}\\
\P^{-1}\I_i&\text{if $r_0+r_1<i\le j$,}\\
\I_i&\text{otherwise.}\endcases$$
Set $m_1=r_1+n-j$.
Expanding $f''$ as a Fourier series supported on even
$T\in((\I_i''\I_{\ell}''\J)^{-1})$, we find that for fixed $\L_1$,

$$\align
&\sum_Y f''|\pmatrix I&Y\\&I\endpmatrix 
\pmatrix C^{-1}\\&^tC\endpmatrix(\tau)\\
&\qquad = \sum_T c_{f''}(T) \e\{TC^{-1}\tau\ ^tC^{-1}\}
\sum_Y \e\{TY\}\\
&\qquad =  \sum_T c_{f''}(T) \e\{TC^{-1}\tau\ ^tC^{-1}\}\\
&\qquad\qquad
\sum_{W_0,W_1,W_2,W_3} \e\{T_0W_0\}\ \e\{\mu T_1W_1\}\ 
\e\{T_2W_2\}\ \e\{T_3W_3\}\endalign$$
where 
$T=\pmatrix T_0&T_2&*&T_3\\ ^tT_2&T_1&*&*\\*&*&*&*\\^tT_3&*&*&*\endpmatrix.$
$T_0$ and $W_0$ are symmetric $r_0\times r_0$ matrices with $T_0$ even,
the $i,\ell$-entry of $W_0$ in $\I_i\I_\ell\J\partial^{-1}$.  Thus
the sum on $W_0$ ($T_0$ fixed) is a complete character sum, yielding
$N(\P)^{r_0(r_0+1)}$ if $T_0\equiv0\ (\mod\P)$, and 0 otherwise.  
Similarly, the sums on $W_2, W_3$ are complete character sums.  So
$$\sum_Y \e\{TY\}=\cases N(\P)^{r_0(r_0+m_1+1)}\sum_{W_1} \e\{\mu T_1W_1\}
&\text{if $T\in((\I_i\I_\ell\J)^{-1})$,}\\
0&\text{otherwise.}\endcases$$
With $(B(x_i,x_\ell))=\ ^tC^{-1}TC^{-1}$, take 
$\L=\I_1 x_1\oplus\cdots\oplus\I_nx_n$.  Then the sum on $W_0, W_2, W_3$
is nonzero if and only if $\L^\J$ is even integral.

Let $(y_1\ldots y_n)=(x_1\ldots x_n)C$ and set
$\O=\I''_1y_1\oplus\cdots\oplus\I_n''y_n$.  Then $(B(y_i,y_\ell))=T$
and $c_{f''}(T)N(\I_1''\cdots\I_n'')^k N(\J)^{nk/2}=c_F(\O^\J)$.

Note that when $c_{f|S^{\ell}(\P)}(T)$ is contributing to the $\L^{\J}$th
coefficient of $f|T_j(\P^2)\in\M_k(\G(\I_1,\ldots,\I_n;\J))$, it
gets normalized by $N(\I_1\ldots\I_n)^kN(\J)^{nk/2}$; when it is
determining a coefficient of $f|S^{\ell}(\P)\in\M_k(\G')$,
$\G'\simeq\G(\P^{-\ell}\I_1,\I_2,\ldots,\I_n;\J)$, it gets normalized
by $N(\P)^{-\ell k}N(\I_1\ldots\I_n)^kN(\J)^{nk/2}$.

So the contribution from $f''$ to the $\L^\J$th Fourier coefficient of
$f'|T_j(\P^2)$ is
$$N(\P)^{j(k-n-1)+r_0(r_0+m_1+1)} \chi(\P)^{2r_2+r_1'}
\sum_{\O,\L_1}
N(\P)^{k(2r_2+r_1)} c_F(\O^\J) \sum_{W_1}\e\{\mu T_{\L_1}W_1\}$$
where $\O$ varies subject to $\P\L\subseteq\O\subseteq\P^{-1}\L$,
$\mult_{\{\L:\O\}}(\P^{-1})=r_0$, $\mult_{\{\L:\O\}}(\P)=r_2$,
$\overline\L_1$ a codimension $n-j$ subspace of $\overline\O_1
\simeq \O\cap\L/\P(\O+\L)$, and $W_1$ varies modulo $\P$ with
$i,\ell$-entry in $\I_i\I_\ell\J\partial^{-1}$, $\P\not|\det W_1$.  Here
$T_{\L_1}=(B(x_{r_0+i},x_{r_0+\ell}))$ is $r_1\times r_1$
where $\L_1=\I_{r_0+1}x_{r_0+1}\oplus\cdots\oplus\I_{r_0+r_1}
x_{r_0+r_1}$.  If $T_1'$
is also a matrix associated to $\L_1$ then there is a change of basis
matrix $G$ whose $i,\ell$-entry lies in $\I_i\I_\ell^{-1}$ so that
$^tG T_1' G=T_{\L_1}$.  Thus $\e\{\mu T_1'W_1\}=\e\{\mu T_{\L_1}
(GW_1\ ^tG)\}$; as $W_1$ varies over invertible matrices modulo $\P$,
so does $GW_1\ ^tG$.  Thus the sum on $W_1$ is independent of the choice
of matrix associated to $\L_1$.

We complete the character sum on $W_1$ by replacing $T_j(\P^2)$ with
 $\widetilde T_j(\P^2)$, then apply Proposition 5.1, where we
consider $\sum_W \e\{\a T_1 W\}$ with $W$ varying over all symmetric
$r_1\times r_1$ integral matrices modulo $\P$.

Notice that for $0\le\ell\le j$,
$$S_{\ell+1}(\P)\cdots S_j(\P)T_\ell(\P^2):
\M_k(\G(\I_1',\ldots,\I_n';\J))\to 
\M_k(\G(\I_1,\ldots,\I_n;\J).$$
Also notice that the number of dimension $r_1=m_1-n+j$ lattices
$\overline\L_1$ containing some dimension $m_1-n+\ell$ lattice
$\overline\D$ is the number of ways to extend $\overline\D$ to a
$j$-dimensional subspace of $\overline\O_1$ (where $\O=\O_0\oplus
\O_1\oplus\O_2$, $\L=\P\O_0\oplus\O_1\oplus\P^{-1}\O_2$).
Extending $\overline\D$ is equivalent to choosing a $j-\ell$ dimensional
subspace of an $n-\ell$ space (here $\dim\overline\O_1=m_1$).
So the number of $\overline\L_1$ containing $\overline\D$ is
$\beta(n-\ell,j-\ell)=\beta_\P(n-\ell,j-\ell)$.  

Note that the coefficient of $f'|\widetilde T_j(\P^2)$ associated to
$\L^{\J}=(\I_1x_1\oplus\cdots\oplus\I_nx_n)^{\J}$
carries a normalizing factor of $N(\I_1\cdots\I_n)^k N(\J)^{nk/2}$,
while the coefficient of $f'|S^{(j)}(\P)$ associated to
$\O^{\J}=(\I_1''y_1\oplus\cdots\oplus\I_n''y_n)^{\J}$
carries a factor of 
$N(\P)^{k(r_0-r_2)}N(\I_1\cdots\I_n)^k N(\J)^{nk/2}.$
Hence, contributing to $c(\L^{\J})$ we have
$$N(\P)^{k(r_2-r_0+j)+r_0(r_0+m_1+1)+r_1(r_1+1)/2-j(n+1)}
\chi(\P)^{2r_2-r_0}\sum_{\L_1}c(\O^{\J})$$
where $\L_1^{\J}$ varies over all totally isotropic codimension $n-j$
sublattices of $\L^{\J}\cap\O^{\J}/\P(\L^{\J}+\O^{\J})$.
Summing over all $\O$, $\P\L\subseteq\O\subseteq\P^{-1}\L$, yields (1).

The proof of (2) is quite similar to the proof of (1), but much simpler,
and so we leave it to the reader.
$\square$
\enddemo

\noindent{\bf Remark.}  As discussed above Proposition 2.2, 
$\L_1^{\J}/\P\L_1^{\J}$
is a quadratic space over $\Ok/\P$.  By \S42 of [6] (for results about
quadratic spaces over fields of characteristic 2, see, for example,
\S5 of [9]), $\L_1^{\J}/\P\L_1^{\J}=R\perp W\perp{\Bbb H}^t$ where
$R=\rad\L_1^{\J}/\P\L_1^{\J}$, $W$ is anisotropic, and 
${\Bbb H}\simeq\pmatrix 0&1\\1&0\endpmatrix$ denotes a hyperbolic plane.
With $U=R\perp W$, Lemma 1.6 of [8] 
and Lemma 4.1 of [9] tell us that the number
of $\ell$-dimensional totally isotropic subspaces of 
$\L_1^{\J}/\P\L_1^{\J}$ is
$$\varphi_{\ell}(\L_1^{\J}/\P\L_1^{\J})
=\sum_a q^{(t-a)(\ell-a)} \delta(d+t-\ell+a+1,a)
\beta(t,a) \varphi_{\ell-a}(U)$$
where $q=N(\P)$, $d=\dim U$, 
$\delta(m,r)=\prod_{i=0}^{r-1} (q^{m-i}+1)$,
$\beta(m,r)=\prod_{i=0}^{r-1}(q^{m-i}-1)/(q^{r-i}-1)$, and
$0\le a \le \ell$.
(Note that Lemma 1.6 of [8] is proved for a quadratic space over
$\Z/p\Z$, but the argument is valid over all finite fields.  When the
characteristic is 2, we replace $Q$ by ${1\over2}Q$; we present
a full discussion of this case in \S5 of [9].)
Also, since $U=R\perp W$ with $R$ totally isotropic and $W$ anisotropic,
any totally isotropic subspace of $U$ is a subspace of $R$.
Thus $\varphi_{\ell-a}(U)=\beta(r,\ell-a)$ where $r=\dim R$.
(So $d=\dim U=r+\dim W$, and since $W$ is anisotropic, $\dim W$ is 0, 1, or
2 by 62:1b of [6] for $q$ odd, and by, e.g., Proposition 5.1 of [9]
for $q$ even.)  Hence
$$\alpha_j(\O,\L) = \sum_a q^{(t-a)(\ell-a)} \delta(d+t+a+1,a) \beta(t,a)
\beta(r,\ell-a)$$
where $0\le a\le\ell$, $\ell=r_1-n+j$.

\head{6.  Lemma on symmetric coprime pairs.}\endhead

\noindent{\bf Definition.}  We say a pair of matrices $(C,D)$ is a symmetric
coprime lower pair for $\G(\I_1,\ldots,\I_n;\J)$ if:
\itemitem{(a)}  $C\ ^tD$ is symmetric;
\itemitem{(b)}  $C\in\J^{-1}\partial\left(\I_i^{-1}\I_j^{-1}\right)$,
$D\in\left(\I_i^{-1}\I_j\right)$;
\itemitem{(c)}  for all prime ideals $\P$,
$$\rank_{\P}\pmatrix \lambda_1\\& \ddots \\ && \lambda_n\endpmatrix
(C,D)\pmatrix 
\matrix \a\lambda_1\\&\ddots\\&&\a\lambda_n\endmatrix \\
& \matrix \a\lambda^{-1}_1\\&\ddots\\&&\a\lambda^{-1}_n\endmatrix \\
\endpmatrix = n$$
where $\a\Op=\J\partial^{-1}\Op$, $\lambda_i\Op=\I_i\Op$.
\smallskip

Since $\G(\I_1,\ldots,\I_n;\J)\approx \G(\Ok,\ldots,\Ok,\I;\J)$ when
$\I\in\cls(\I_1\cdots\I_n)$, the following technical lemma focuses on
$\G(\Ok,\ldots,\Ok,\I;\J)$, but the conclusion (c) is valid for 
$\G(\I_1,\ldots,\I_n;\J)$.

Corresponding definitions and results hold for symmetric coprime upper,
right-hand, and left-hand pairs for $\G(\I_1,\ldots,\I_n;\J)$.

\proclaim{Lemma 6.1} Let $\I,\J$ be fractional ideals; set
$$\Cal G=\{ E\in\pmatrix I\\&\I\endpmatrix 
\Ok^{n,n}\pmatrix I\\&\I^{-1}\endpmatrix :\ 
\det E=1\ \}.$$
Say $C\in\J^{-1}\partial\pmatrix I\\&\I^{-1}\endpmatrix \Ok^{n,n}
\pmatrix I\\&\I^{-1}\endpmatrix$,
$D\in
\pmatrix I\\&\I^{-1}\endpmatrix \Ok^{n,n}
\pmatrix I\\&\I\endpmatrix$ are a symmetric coprime lower pair for
$\G(\Ok,\ldots,\Ok,\I;\J)$.
\item{(a)} Let $\Q$ be a prime ideal, and choose $\lambda\in\I$,
$\alpha\in\J\partial^{-1}$ such that $\ord_{\Q}\lambda=\ord_{\Q}\I$,
$\ord_{\Q}\alpha=\ord_{\Q}\J\partial^{-1}.$  Then 
 there is some $E\in\Cal G$ so that
$\a\pmatrix I\\&\lambda\endpmatrix CE \pmatrix I\\&\lambda\endpmatrix
\equiv (C_0,0)\ (\mod \Q)$
where $C_0$ is $n\times r$ with $\rank_{\Q}C_0=r$.
\item{(b)}  There is a matrix $M\in\G(\Ok,\ldots,\Ok,\I;\J)$ so that
$(C|D)M=(C'|D')$ with $\det C',\det D'\not=0$ and
$\left( (\J\partial^{-1})^n\I^2\det C',\det D'\right) = 1.$
\item{(c)}  There are $n\times n$ matrices $A,B$ so that
$\pmatrix A&B\\ C&D\endpmatrix \in \G(\Ok,\ldots,\Ok,\I;\J).$
\endproclaim

\demo{Proof}  It suffices to establish the claims for 
$\G(\Ok,\ldots,\Ok,\P;\J)$ where $\P$ is a prime ideal in 
$\cls\I$.  Choose $\a\in\J\partial^{-1}$ so that
$\ord_{\Q}\a=\ord_{\Q}\J\partial^{-1}$ whenever $\ord_{\Q}\J\partial^{-1}
\not=0$.  Choose $\lambda\in\P - \P^2$, $\mu\in\P^{-1} - \Ok$ so that
$\lambda\mu\equiv 1\ (\mod\P)$.

\noindent
(a)  Let $\widetilde C=\a\pmatrix I\\&\lambda\endpmatrix C
\pmatrix I\\&\lambda\endpmatrix=(\ul c_1\ldots \ul c_n)$ (so $\ul c_i$ is the 
$i$th
column of $\widetilde C$).  Then there is some $E'\in GL_{n-1}(\Ok)$ so that
$$(\ul c_1\ldots \ul c_{n-1})E\equiv (C_0',0)\ (\mod\P)$$
where $C_0'$ is $n\times r'$, $r'=\rank_{\P}(\ul c_1\ldots \ul c_{n-1})$.
So $E_1=\pmatrix E'\\&1\endpmatrix \in \Cal G$, and
$$\a\pmatrix I\\&\lambda\endpmatrix CE_1\pmatrix I\\&\lambda\endpmatrix
=\widetilde CE_1\equiv (C_0',0,\ul c_n)\ (\mod\P).$$

First suppose $\ul c_n\in\spn_{\P}C_0'.$  Thus there exist
$\g_1,\ldots,\g_r\in\Ok$ so that $C_0'\underline\g\equiv -\ul c_n\ (\mod\P)$
where $\underline\g=\pmatrix \g_1\\ \vdots\\ \g_r\endpmatrix.$  Let
$E=E_1\pmatrix I&&\mu\underline\g\\&I\\&&1\endpmatrix.$
Then
$$\align
\a\pmatrix I\\&\lambda \endpmatrix CE\pmatrix I\\&\lambda\endpmatrix
&=\widetilde CE_1\pmatrix I\\&(\mu\lambda)^{-1}\endpmatrix
\pmatrix I&&\underline\g\\&I\\&&1\endpmatrix
\pmatrix I\\&\mu\lambda\endpmatrix \\
&\equiv\widetilde CE_1\pmatrix I\\&(\mu\lambda)^{-1}\endpmatrix
\pmatrix I&&\underline\g\\&I\\&&1\endpmatrix
\pmatrix I\\&\mu\lambda\endpmatrix\ (\mod\P) \\
&\equiv (C_0',0)\ (\mod\P),\endalign$$
proving (a) in the case $\ul c_n\in\spn_{\P}C_0'$.

So suppose $\ul c_n\not\in \spn_{\P}C_0'$.  If $r'=n-1$ then we are done,
as then $\rank_{\P}\widetilde C=n$ and we can take $E=I$.  So let us also
suppose $r'<n-1$.  Choose $\delta\in\P$ so that for all primes $\Q|\mu\lambda$,
$\Q\not|\delta$.  Thus $(\mu\lambda,\delta)=1$, so there are $u,v\in\Ok$ so that
$v\delta-u\mu\lambda=1.$  Thus
$E_2=\pmatrix I\\&\matrix \nu&\mu\eta\\ \lambda&\delta\endmatrix \endpmatrix
\in\Cal G,$ and
$$\a\pmatrix I\\&\lambda\endpmatrix CE_1E_2\pmatrix I\\&\lambda\endpmatrix
\equiv \widetilde C E_1\cdot \pmatrix I\\&\mu\endpmatrix E_2
\pmatrix I\\&\lambda\endpmatrix\ (\mod\P)$$
where $\ul c'_i=\ul c_i$ for $i<n-1$, $\ul c'_{n-1}\equiv \ul c_n\ (\mod\P)$ and
$\ul c'_n\equiv 0\ (\mod\P).$

Let $E_3=\pmatrix *\\&1\endpmatrix$ be a permutation matrix that interchanges
columns $r'$ and $n-1$.  So $E=E_1E_2E_3\in\Cal G$, and (a) is proved.

\noindent (b)  
Let $\Cal A = (\J\partial^{-1})^n\cdot \I^2\cdot \det C
\subseteq\Ok.$
Take $\lambda\in\I$, $\eta\in\I^{-1}$, $\alpha\in\J\partial^{-1}$,
$\mu\in\J^{-1}\partial$ so that 
$\lambda\eta\equiv\alpha\nu\equiv 1\ (\mod \Q)$.

{\bf Case 1:}  Say $\det D\not=0$, $\det D\not\in\Ok^{\times}$.
Part (a) says there exists some $E\in\Cal G$ so that
$$\alpha\pmatrix I\\&\lambda\endpmatrix CE \pmatrix I\\&\lambda\endpmatrix
\equiv (C_0,0)\ (\mod \Q),$$
with $C_0$ $n\times r$, $\rank_{\Q}C_0=r$.  Thus
$\pmatrix E\\&^tE^{-1}\endpmatrix\in\G(\Ok,\ldots,\Ok,\I;\J)$,
and
$$\pmatrix I\\&\lambda\endpmatrix (C|D) \pmatrix E\\&^tE^{-1}\endpmatrix
\pmatrix \alpha I\\&\alpha\lambda\\&&I\\&&&\eta\endpmatrix
\equiv (C_0,0)D_0,D_1)\ (\mod \Q);$$
here $D_0\subseteq\spn_{\Q}C_0$ (since $C\ ^tD$ is symmetric; see
the proof of Lemma 7.2 in [5]), and thus $\rank_{\Q}(C_0|D_1)=n$.

Set $W=\mu\pmatrix I\\&\eta\endpmatrix \pmatrix 0_r\\&I_{n-r}\endpmatrix
\pmatrix I\\&\eta\endpmatrix$; hence
$$M_{\Q}=\pmatrix E\\&^tE^{-1}\endpmatrix \pmatrix I&0\\W&I\endpmatrix
\in\G(\Ok,\ldots,\Ok,\I;\J).$$
Also, with $(C'|D')=(C|D)M_{\Q}$, $\det D''=\det D$ and
$$\align
& \pmatrix I\\&\lambda\endpmatrix (C'|D') \pmatrix \alpha I\\&\alpha\lambda\\
&&I\\&&&\eta\endpmatrix\\
&\qquad
\equiv \pmatrix I\\&\lambda\endpmatrix (C|D) \pmatrix E\\&^tE^{-1}\endpmatrix
\pmatrix \alpha I\\&\alpha\lambda\\&&I\\&&&\eta\endpmatrix
\pmatrix I\\&I\\0_r&&I\\&I&&I\endpmatrix \\
&\qquad\qquad (\mod \Q).
\endalign$$
Thus $\rank_{\Q} \alpha\pmatrix I\\&\lambda\endpmatrix C'
\pmatrix I\\&\lambda\endpmatrix =n$, and hence $\Q\not|(\J\partial^{-1})^n
\I^2\det C'$.

Only finitely many prime ideals $\Q$ divide $\det D$, so repeating
this process for all such $\Q$ yields a pair $(C'|D')=(C|D)M$,
$M\in\G(\Ok,\ldots,\Ok,\I;\J)$, $\det C',\det D'\not=0$,
and $\big( (\J\partial^{-1})^n\I^2\det C',\det D')=1.$

{\bf Case 2:}  Say $\det D\in\Ok^{\times}$.  If $\det C\not=0$ then we are
done; so suppose $\det C=0$.  Then let $\Q$ be any prime ideal; following
the algorithm in Case 1, we produce $M_{\Q}\in\G(\Ok,\ldots,\Ok,\I;\J)$ so
that with $(C'|D')=(C|D)M_{\Q}$, we have $\det D'=\det D$ and 
$\Q\not|\det C'$ (so $\det C'\not=0$).  This completes this case.

{\bf Case 3:}  Say $\det D=0$.  Choose a prime $\Q$; with $\lambda,\eta,
\alpha,\mu$ as in Case 1, we know there is some $E\in\Cal G$ so that
$$\pmatrix I\\&\lambda\endpmatrix (C|D) \pmatrix E\\&^tE^{-1}\endpmatrix
\pmatrix \alpha I\\&\alpha\lambda\\&&I\\&&&\eta\endpmatrix
\equiv (C_0,0|D_0,D_1)\ (\mod \Q)$$
with $\rank_{\Q}(C_0|D_1)=n$, $D_0\subseteq\spn_{\Q}C_0.$
Thus for some 
$Y\in\Ok^{r,r}$, $\rank_{\Q}(C_0Y+D_0,D_1)=n$.  Note that $Y$ is uniquely 
determined modulo $\Q$, and that $C_0\ ^tD_0$ is symmetric modulo $\Q$.
Since $C_0Y\ ^tC_0\equiv -D_0\ ^tC_0\ (\mod \Q)$, we can choose $Y$
to be symmetric in $\Ok^{r,r}$.

Set
$$Y'=\alpha\pmatrix I\\&\lambda\endpmatrix \pmatrix Y\\&0\endpmatrix
\pmatrix I\\&\lambda\endpmatrix;$$
hence $M=\pmatrix E\\&^tE^{-1}\endpmatrix \pmatrix I&Y'\\&I\endpmatrix
\in\G(\Ok,\ldots,\Ok,\I;\J).$
Thus with $(C'|D')=(C|D)M$, we have
$$\pmatrix I\\&\lambda\endpmatrix(C'|D')\pmatrix\alpha I\\&\alpha\lambda\\
&&I\\&&&\eta\endpmatrix
\equiv (C_0,0|C_0Y+D_0,D_1)\ (\mod\Q).$$
Hence $\rank_{\Q}\pmatrix I\\&\lambda\endpmatrix
 D'\pmatrix I\\&\eta\endpmatrix=n$,
so $\Q\not|\det D'$ (and thus $\det D'\not=0$).  This reduces this situation
to one of the previous cases.

\noindent (c) Set $\kappa=\det D$; choose $\lambda\in\Cal A$ so that
$(\kappa,\lambda)=1$, and $\eta\in\Ok$ so that $\eta\equiv \kappa^{-1}
\ (\mod \lambda)$ (this is possible by the Chinese Remainder Theorem).
Since $C\in\J^{-1}\partial \pmatrix I\\&\I^{-1}\endpmatrix \Ok^{n,n}
\pmatrix I\\&\I^{-1}\endpmatrix,$ we have
$\lambda\ ^tC^{-1}\in\J\partial^{-1}\pmatrix I\\&\I\endpmatrix \Ok^{n,n}
\pmatrix I\\&\I\endpmatrix.$
Also, $\lambda|(\eta\kappa-1)$, so
$$B=(\eta\kappa-1)\ ^tC^{-1} \in \J\partial^{-1}
\pmatrix I\\&\I\endpmatrix \Ok^{n,n}
\pmatrix I\\&\I\endpmatrix.$$
(Note that locally everywhere:
 $C\in\J^{-1}\partial \pmatrix I\\&\I^{-1}\endpmatrix C_0
\pmatrix I\\&\I^{-1}\endpmatrix,$ $C_0\in\Ok^{n,n}$,
$\det C_0=(\J\partial^{-1})^n\I^2\det C$; hence
 $C^{-1}\in {1\over \det C_0}
\J\partial^{-1} \pmatrix I\\&\I\endpmatrix \Ok^{n,n}
\pmatrix I\\&\I\endpmatrix.$)

Set $A=\eta\kappa\ ^tD^{-1}.$  We have $\kappa=\det D$, $\eta\in\Ok$,
so $\eta\kappa\ ^tD^{-1}\in\pmatrix I\\&\I\endpmatrix\Ok^{n,n}
\pmatrix I\\&\I^{-1}\endpmatrix.$  Thus $\pmatrix A&B\\C&D\endpmatrix$
is a candidate for $\G(\Ok,\ldots,\Ok,\I;\J)$.  To see this matrix 
indeed lies in this group, first note that $A\ ^tD-B\ ^tC=I$, and
by assumption, $C\ ^tD$ is symmetric.  Finally, substituting for
$A$ and $B$, we get $A\ ^tB=\eta\kappa(\eta\kappa-1)\ ^tD^{-1}C^{-1}$;
since $^tD^{-1}C^{-1}=(C\ ^tD)^{-1}$ is symmetric, so is $A\ ^tB$.
$\square$
\enddemo

\enddocument